\documentclass[a4paper,11pt]{amsart}
\usepackage[OT2,OT1]{fontenc}
\usepackage[latin1]{inputenc}
\usepackage{amsmath}
\usepackage{amsthm}
\usepackage{amscd}
\usepackage{amssymb}
\usepackage{pb-diagram, pb-xy}
\usepackage[all]{xy}
\usepackage{hyperref}
\hypersetup{colorlinks=true}

\theoremstyle{plain}
	\newtheorem{mainthm}{Theorem}
	\newtheorem{thm}{Theorem}[section]
	\newtheorem{cor}[thm]{Corollary}
	\newtheorem{lem}[thm]{Lemma}
	\newtheorem{prop}[thm]{Proposition}
\theoremstyle{definition}
	\newtheorem{defn}[thm]{Definition}
	
	\newtheorem{que}[thm]{Question}
	\newtheorem{para}[thm]{--}

\numberwithin{equation}{section}

\newcommand\cyr{%
\renewcommand\rmdefault{wncyr}%
\renewcommand\sfdefault{wncyss}%
\renewcommand\encodingdefault{OT2}%
\normalfont\selectfont}
\DeclareTextFontCommand{\textcyr}{\cyr}

\newcommand{\rnode}[1]{\node{\raisebox{-7mm}{$\displaystyle #1$}}}

\newcommand{\angl}[1]{\left\langle #1\right\rangle}

\newcommand{\eps}{\varepsilon}
\newcommand{\tq}{\: | \:}
\newcommand{\pt}{\mbox{ for all }}

\newcommand{\qqet}{\qquad\mbox{and}\qquad}

\newcommand{\Scha}{\textcyr{X}}

\DeclareMathOperator{\colim}{colim}%
\DeclareMathOperator{\End}{End}
\DeclareMathOperator{\Gal}{Gal}
\DeclareMathOperator{\GL}{GL}
\DeclareMathOperator{\Hom}{Hom}
\DeclareMathOperator{\id}{id} %
\DeclareMathOperator{\im}{im} %
\DeclareMathOperator{\Pic}{Pic} %
\DeclareMathOperator{\red}{red}
\DeclareMathOperator{\res}{res}
\DeclareMathOperator{\spec}{spec}
\DeclareMathOperator{\tr}{tr}

\newcommand{\hotimes}{\:\widehat\otimes\:}
\newcommand{\Iell}{\textup{B}_\ell}
\newcommand{\Tell}{\textup{T}_\ell}

\newcommand{\Vell}{\textup{V}_{\!\ell}}
\newcommand{\Vnot}{\textup{V}_{\!0}}

\newcommand{\upet}{\textup{\'et}}
\newcommand{\uptor}{\textup{tor}}

\newcommand{\IC}{\mathbb{C}}
\newcommand{\IF}{\mathbb{F}}
\newcommand{\IG}{\mathbb{G}}
\newcommand{\IL}{\mathbb{L}}

\newcommand{\IQ}{\mathbb{Q}}
\newcommand{\IR}{\mathbb{R}}
\newcommand{\IZ}{\mathbb{Z}}
\newcommand{\QZ}{\IQ/\IZ}

\newcommand{\cC}{\mathcal C}
\newcommand{\cD}{\mathcal D}
\newcommand{\cF}{\mathcal F}
\newcommand{\cM}{\mathcal M}
\newcommand{\cO}{\mathcal O}

\newcommand{\cAb}{\mathcal Ab}
\newcommand{\cHom}{\mathcal Hom}
\newcommand{\cExt}{\mathcal Ext}

\newcommand{\fc}{\mathfrak c}
\newcommand{\fd}{\mathfrak d}
\newcommand{\fl}{\mathfrak l}

\newcommand{\fn}{\mathfrak n}
\newcommand{\fp}{\mathfrak p}

\renewcommand{\to}{\xrightarrow{\quad}}
\renewcommand{\mapsto}{\longmapsto}

\title{On the second Tate--Shafarevich group of a 1--motive}
\author{Peter Jossen}
\date{\today}

\begin{document}

\begin{abstract}
We prove finiteness results for Tate--Shafarevich groups in degree 2 associated with 1--motives. We give a number theoretic interpretation of these groups, rely them to Leopoldt's conjecture, and present an example of a semiabelian variety with an infinite Tate--Shafarevich group in degree 2. We also establish an arithmetic duality theorem for 1--motives over number fields which complements earlier results of Harari and Szamuely.
\end{abstract}

\maketitle

\section*{Introduction and overview}

\begin{par}
Let $k$ be a number field and let $X$ be a commutative group scheme over $k$. The Tate--Shafarevich group $\Scha^i(k,X)$ of $X$ is the subgroup of the \'etale cohomology group $H^i(k,X)$ consisting of those elements which restrict to zero over each completion of $k$. These groups are among the most fundamental invariants associated with commutative group schemes over number fields, and their vanishing is by definition the obstruction to various local--to--global principles. 
\end{par}

\begin{par}
If the group scheme $X$ is given by a finitely generated discrete group with continuous Galois action, or if $X$ is a group of multiplicative type, then $\Scha^i(k,X)$ is finite for all $i$ (\cite{MilneADT} Theorem I.4.20 and \cite{Neukirch00}, Theorem 8.6.8). It is widely conjectured that if $A$ is an abelian variety over $k$, then the group $\Scha^1(k,A)$ is finite, and it is known that for $i\neq 1$ the group $\Scha^i(k,A)$ is trivial. This is a nontrivial statement for $i=2$, indeed, the proof of Corollary I.6.24 in \cite{MilneADT} shows that the vanishing of $\Scha^2(k,A)$ is essentially equivalent to the positive answer to the congruence subgroup problem for the abelian variety dual to $A$, given by Serre in \cite{Serre1} and \cite{Serre2}. 
\end{par}

\begin{par}
An evident generalisation of these finiteness results would be to show that $\Scha^i(k,G)$ is finite for semiabelian varieties $G$ over $k$, i.e.\ when $G$ is an extension of an abelian variety $A$ by a torus. A simple d\'evissage shows that $\Scha^1(k,G)$ is finite (\cite{Hara05}, Lemma 4.11), provided $\Scha^1(k,A)$ is finite. The situation is more complicated for $i=2$, and surprisingly it turns out that the groups $\Scha^2(k,G)$ are not always finite.
\end{par}

\vspace{2mm}
\begin{mainthm}\label{Thm:Main1}
There exists a semiabelian variety $G$ over $\IQ$ such that the group $\Scha^2(\IQ,G)$ contains a subgroup isomorphic to $\QZ$, and in particular is infinite.
\end{mainthm}

\vspace{2mm}
\begin{par}
A 1--motive $M$ over a number field $k$ is a two term complex of group schemes $M=[Y\to G]$ over $k$ placed in degrees $-1$ and $0$, where $Y$ is given by a finitely generated free discrete group with continuous Galois action, and where $G$ is a semiabelian variety. It was asked in \cite{Hara05}, Remark 4.13 whether for all 1--motives $M$ the group $\Scha^2(k,M)$ is finite. By Theorem 1 we already know that this is not always the case even for 1--motives of the form $[0\to G]$ over $\IQ$. Our second result shows that even for very simple 1--motives it might be difficult to decide whether $\Scha^2(k,M)$ is finite or not (assuming the conservation law of difficulity).
\end{par}

\vspace{2mm}
\begin{mainthm}\label{Thm:Main2}
If the group $\Scha^2(k,M)$ is finite for all 1--motives of the form $M= [\IZ^r \to \IG_m^s]$ over $k$, then Leopoldt's conjecture holds for $k$ (and all prime numbers).
\end{mainthm}

\vspace{2mm}
\begin{par}
Our third result provides conditions on a 1--motive which ensure that $\Scha^2(k,M)$ is finite. It is most conveniently expressed as a duality theorem. Classical global arithmetic duality theorems are statements about the existence and nondegeneracy of canonical pairings between certain Tate--Shafarevich groups. Let $X$ be a group of multiplicative type over $k$, and denote by $X^\vee$ its group of characters. The Poitou--Tate Duality Theorem states that there is a natural, perfect pairing of finite groups
$$\Scha^i(k,X) \times \Scha^{3-i}(k,X^\vee) \to \QZ$$
(\cite{MilneADT} Theorem I.4.20 and \cite{Neukirch00}, Theorem 8.6.8). The analogue of this duality theorem for abelian varieties is the Cassels--Tate Duality Theorem. It states that for an abelian variety $A$ over $k$ with dual $A^\vee$, there is a canonical pairing 
$$\Scha^i(k,A) \times \Scha^{2-i}(k,A^\vee) \to \QZ$$
whose left and right kernels are the maximal divisible subgroups (\cite{MilneADT} Theorem I.6.26). Conjecturally it is a perfect pairing of finite groups. 
\end{par}

\begin{par}
The idea to unify and generalise these arithmetic duality theorems to duality theorems for 1--motives is due to Harari and Szamuely and appeared first in \cite{Hara05}. Deligne constructed for each 1--motive $M$ a dual 1--motive $M^\vee$. In \cite{Hara05} the authors show that for a 1--motive $M$ over a number field $k$ there is a canonical pairing
$$\Scha^1(k,M) \times \Scha^1(k,M^\vee) \to \QZ$$
which is nondegenerate modulo divisible subgroups and generalises the Cassels--Tate pairing. Harari and Szamuely also construct a pairing between a certain modification of $\Scha^0(k,M)$ and $\Scha^2(k,M^\vee)$ and show that it is nondegenerate modulo divisible subgroups. However, this modified $\Scha^0(k,M)$ remains somehow uncontrollable, and the resulting generalised pairing does not seem to be very useful (the statement of Theorem 0.2 in \cite{Hara05} was rectified in \cite{Hara05E}).
\end{par}

\vspace{2mm}
\begin{mainthm}\label{Thm:Main3}
Let $k$ be a number field and let $M = [u: Y\to G]$ be a 1--motive over $k$ with dual $M^\vee$. There exists a natural pairing
\begin{equation*}
\Scha^0(k,M) \times \Scha^2(k,M^\vee) \to \QZ \tag{$\ast$}
\end{equation*}
generalising the Poitou--Tate pairing for finitely generated Galois modules and tori. The group $\Scha^0(k,M)$ is finite and the pairing $(\ast)$ is nondegenerate on the left. If the semiabelian variety $G$ is a an abelian variety or a torus such that the $\IQ$--algebra $\End_{\overline k}(G)\otimes\IQ$ is a product of division algebras, then the pairing $(\ast)$ is a perfect pairing of finite groups.
\end{mainthm}

\vspace{2mm}
\begin{par}
It was already shown in \cite{Hara05} that $\Scha^0(k,M)$ is finite. The finiteness results stated in the second part of the theorem are new, and are also the essential part of the theorem. Equivalently, the condition on $G$ is that over an algebraic closure of $k$ either $G$ is the multiplicative group, or an abelian variety isogenous to a product of pairwise nonisogenous simple abelian varieties. Our proof uses techniques developed by Serre in his work on the congruence subgroup problem for abelian varieties (\cite{Serre1}, \cite{Serre2}).
\end{par}

\vspace{4mm}
\begin{par}{\bf Overview:}
In section \ref{Sec:About1Mot} we rehearse 1--motives and $\ell$--adic realisations, which will play a prominent role throughout this paper. In section \ref{Sec:PairingScha0Scha2} we construct a duality pairing which relates the $\ell$--adic realisation of a 1--motive with the second Tate--Shafarevich group of its dual, and obtain the pairing $(\ast)$ of Theorem \ref{Thm:Main3}. In \ref{Sec:H1ast} we compute the cohomology of some $\ell$--adic Lie groups associated with 1--motives, and in \ref{Sec:Finiteness} we use these computations to prove the finiteness statements in Theorem \ref{Thm:Main3}. We conclude the proof of Theorem \ref{Thm:Main3} in section \ref{Sec:TorsionOfSchaTellM}. In \ref{Sec:Leopoldt} and \ref{Sec:InfiniteScha2} we prove the Theorems \ref{Thm:Main2} and \ref{Thm:Main1} respectively. There remain several open questions and unsolved problems, which I state in the last section.
\end{par}

\vspace{4mm}
\begin{par}{\bf Acknowledgments:}
Some parts of this work are taken from my Ph.D. thesis under the direction of Tam\'as Szamuely. I wish to thank him for suggesting the problems and for patiently guiding me. I am also much indebted to David Harari and Gergely Harcos for numerous useful comments. While completing this work I was generously supported by the Central European University and the R\'enyi Institute of Mathematics in Budapest, the University of Regensburg and the Fondation Math\'ematique Jacques Hadamard (FMJH). 
\end{par}

\tableofcontents

\vspace{4mm}
\section{About 1--motives and their realisations}\label{Sec:About1Mot} 

\begin{par}
In this section we rehearse the relevant facts about classical 1--motives and their realisations defined by Deligne in \cite{DeligneHodgeIII}, \S 10.
\end{par}

\vspace{4mm}
\begin{para}
Throughout this section, $S$ is a noetherian regular scheme, $\cF_S$ stands for the category of sheaves of commutative groups on the small fppf site over $S$ and $\cD\cF_S$ for the derived category of $\cF_S$. We identify commutative group schemes over $S$ with objects of $\cF_S$ via the functor of points. In particular, we say that an fppf sheaf on $S$ is an abelian scheme, a torus or a finite flat group scheme if it can be represented by such. By a \emph{lattice over $S$} we mean an object of $\cF_S$ which is locally isomorphic to a finitely generated free $\IZ$--module. Notice that if $S$ is the spectrum of a field, then $Y$ may be regarded as a finitely generated group on which the absolute Galois group acts continuously. 
\end{para}

\vspace{4mm}
\begin{defn}\label{Def:1MotifFleches}
A \emph{$1$--motive} over $S$ is a diagram
$$M = \left[\begin{diagram} \setlength{\dgARROWLENGTH}{4mm}
\node[3]{Y}\arrow{s,l}{u}\\
\node{0}\arrow{e}\node{T}\arrow{e}\node{G}\arrow{e}\node{A}\arrow{e}\node{0}
\end{diagram}\right]$$
in the category $\cF_S$, where $Y$ is a lattice, $T$ a torus, and $A$ an abelian scheme. A morphism of $1$--motives $\varphi:M_1 \to M_2$ is a morphism between diagrams. The \emph{complex associated with $M$} is the complex of fppf--sheaves $[M] := [Y \to G]$, placed in degrees $-1$ and $0$. We denote by $\cM_{1,S}$ the category of 1--motives over $S$.
\end{defn}

\vspace{4mm}
\begin{para}
Observe that the sheaf $G$ is representable. Indeed, we may look at it as an $T$--torsor over $A$, and since $T$ is affine, representability follows from \cite{MilneEC}, Theorem III.4.3a. Later on, 1--motives will often be given by their associated complexes, and morphisms accordingly by commutative squares. This is also custom in the literature, and justified by the fact that there are no nontrivial morphisms from a torus to an abelian scheme.
\end{para}

\vspace{4mm}
\begin{para}\label{Par:WeightFiltration}
We say that a sequence of morphisms of 1--motives  $0 \to M_1 \to M_2 \to M_3 \to 0$ is a \emph{short exact sequence} if the induced sequences of lattices, tori and abelian schemes are exact in $\cF_S$. Such a short exact sequence of 1--motives yields then an exact triangle
$$[M_1]\to[M_2]\to[M_3] \to [M_1][1]$$
in the derived category $\cD\cF_S$. With a 1--motive $M$ over $S$ are naturally associated several short exact sequences coming from the \emph{weight filtration} on $M$. This is the natural three term filtration given by $W_i M = 0$ if $i \leq -3$ and $W_i M  = M$ if $i \geq 0$ and
$$W_{-2}M := \left[\begin{diagram}
\setlength{\dgARROWLENGTH}{4mm}%
\node[3]{0}\arrow{s}\\
\node{\!0\!}\arrow{e}\node{\!T\!}\arrow{e,=}\node{\!T\!}\arrow{e}\node{\!0\!}
\arrow{e}\node{\!0\!}
\end{diagram}\right] \qqet
W_{-1}M := \left[\begin{diagram}
\setlength{\dgARROWLENGTH}{4mm}%
\node[3]{0}\arrow{s}\\
\node{\!0\!}\arrow{e}\node{\!T\!}\arrow{e}\node{\!G\!}\arrow{e}\node{\!A\!}
\arrow{e}\node{\!0\!}
\end{diagram}\right]
$$
Although 1--motives do not form an abelian category, the quotients $M/W_iM$ make sense in the obvious way.
\end{para}

\vspace{4mm}
\begin{defn}
Let $M$ be a 1--motive over $S$ and let $\ell$ be a prime number which is invertible on $S$. The \emph{$\ell$--adic Tate module} and the \emph{$\ell$--divisible Barsotti--Tate group} associated with $M$ are the smooth $\ell$--adic sheaf, respectively $\ell$-divisible group 
$$\Tell(M) := \lim_{i\geq 0}H^1([M]\otimes^\IL\IZ/\ell^i\IZ) \qquad \qquad \Iell(M) := \underset{i\geq 0}\colim\:H^1([M]\otimes^\IL\IZ/\ell^i\IZ)$$
over $S$, where the derived tensor product is taken in the derived category $\cD\cF_S$.
\end{defn}

\vspace{4mm}
\begin{para}
By construction $\Tell M$ only depends on the complex $[M] = [Y\xrightarrow{\:\:u\:\:} G]$ up to quasi--isomorphism, and the assignment $M \mapsto \Tell M$ is functorial. Using the flat resolution $\IZ\xrightarrow{\:\:\ell^i\:\:}\IZ$ of $\IZ/\ell^i\IZ$, we see that the object $[M]\otimes^\IL\IZ/\ell^i\IZ$ of $\cD\cF_S$ is given by the bounded complex
$$\cdots \to 0 \to Y \xrightarrow{\:y\mapsto (\ell^iy,u(y))\:} Y\oplus G \xrightarrow{\:(y,g)\mapsto u(y)-\ell^ig\:} G \to 0 \to \cdots$$
supported in degrees $0,1$ and $2$. For $n\neq 1$ we have $H^n([M]\otimes^\IL\IZ/\ell^i\IZ)=0$, because $Y$ is torsion free and $G$ is divisible as a sheaf. Hence the object $[M]\otimes^\IL\IZ/\ell^i\IZ$ of $\cD\cF_S$ is homologically concentrated in degree $1$. Given a category $\cC$ and a functor $F: \cD\cF_S \to \cC$, we define 
$$F(\Tell M) := \lim_{i\geq 0} F([M]\otimes^\IL\IZ/\ell^i\IZ[-1] ) \qqet F(\Iell M) := \underset{i\geq 0}\colim\:  F([M]\otimes^\IL\IZ/\ell^i\IZ[-1])$$
Depending on the context, these are either limit systems in $\cC$ or actual objects of $\cC$, provided limits and colimits exist in $\cC$.
\end{para}

\vspace{4mm}
\begin{para}\label{Par:TellMExplicit}
Suppose $S$ is connected, and let $\spec(\overline k)=\overline s \to S$ be a geometric point where $\overline k$ is an algebraic closure of the residue field $k$ at the scheme point underlying $\overline s$. We can describe the finite, locally constant group schemes $H^1([M]\otimes^\IL\IZ/\ell^i\IZ)$ in terms of finite $\pi_1 := \pi_1^\upet(\overline s,S)$--modules as follows: The underlying group is given by
$$\frac{\{(y,P) \in Y(\overline k) \times G(\overline k)\tq \ell^iP=u(y)\}}{\{(\ell^iy,u(y))\tq y\in Y(\overline k)\}}$$
and the action of $\pi_1$ is induced by the action of $\pi_1$ on $\overline k$. Taking the limit over $i\geq 0$, we find the description of the $\ell$--adic sheaf $\Tell M$ as a $\pi_1$--module. The short exact sequence of 1--motives coming from the weight filtration $0\to [0\to G] \to M \to [Y\to 0] \to 0$ induces a sequence of $\ell$--adic sheaves, respectively continuous $\pi_1$--representations
$$0\to \Tell G \to \Tell M \to Y \otimes \IZ_\ell \to 0$$
which is exact because $G(\overline k)$ is a $\ell$--divisible group. Observe that, given $y\in Y$, a preimage of $y\otimes 1$ in $\Tell M$ is given by a sequence $(y,P_i)_{i\geq 0}$ with $P_0 = u(y)$ and $\ell P_i = P_{i-1}$ for $i\geq 1$.
\end{para}

\vspace{4mm}
\begin{para}\label{Par:GroupCompletions}
Let $A$ be a commutative group. We consider the following four operations on $A$ relative to the prime $\ell$
$$\begin{array}{cc}
A\hotimes\IZ_\ell := \displaystyle\lim_{i\geq 0}A/\ell^i A \qquad&\qquad \Tell A := \displaystyle\lim_{i\geq 0}A[\ell^i]\\[3mm]   
A[\ell^\infty] := \underset{i\geq 0}{\colim}\: A[\ell^i] \qquad&\qquad A \otimes \IQ_\ell/\IZ_\ell := \underset{i\geq 0}{\colim}\:A/\ell^i A
\end{array}$$
These are the $\ell$--adic completion, the $\ell$--adic Tate module, extraction of $\ell$--torsion and tensorisation with $\IQ_\ell/\IZ_\ell$. These four operations are related, as follows. Given a short exact sequence of commutative groups $0\to A \to B \to C \to 0$, there is a long exact sequence of $\IZ_\ell$--modules
$$0 \to \Tell A \to \Tell B \to \Tell C \to A\hotimes\IZ_\ell \to B\hotimes\IZ_\ell \to C\hotimes\IZ_\ell \to 0$$
coming from the snake lemma, identifying $-\hotimes \IZ_\ell$ with the first right derived functor of the Tate module functor $\Tell(-)$ and vice versa. Similarly, there is a six term exact sequence of $\ell$--torsion groups
$$0 \to A[\ell^\infty] \to B[\ell^\infty] \to C[\ell^\infty] \to A\otimes\IQ_\ell/\IZ_\ell\to B\otimes\IQ_\ell/\IZ_\ell\to
C \otimes\IQ_\ell/\IZ_\ell\to 0$$
identifying $(-)[\ell^\infty]$ with the first left derived functor of $-\otimes\IQ_\ell/\IZ_\ell$ and vice versa. Given a bilinear pairing of commutative groups $A\times B \to \QZ$, these operations induce pairings
$$A\hotimes\IZ_\ell \times B[\ell^\infty] \to \QZ \qqet \Tell A \times (B \otimes\IQ_\ell/\IZ_\ell)\to \QZ$$
If the original pairing was nondegenerate, these are nondegenerate pairings as well. Most of the time we shall deal with commutative groups on which the multiplication--by--$\ell$ has finite kernel and cokernel. For such a group $A$, the $\IZ_\ell$--modules $A \hotimes \IZ_\ell$ and $\Tell A$ are finitely generated, and the torsion groups $A \otimes \IQ_\ell/\IZ_\ell$ and $A[\ell^\infty]$ are of cofinite type (meaning that their Pontryagin duals are finitely generated as $\IZ_\ell$--modules), and there is an isomorphism of finite groups $(A\hotimes \IZ_\ell)[\ell^\infty] \cong A[\ell^\infty] \hotimes\IZ_\ell$. Nondegenerate pairings of such groups induce perfect pairings of topological groups.
\end{para}

\vspace{4mm}
\begin{prop}\label{Pro:FiltrationTate}
Let $F: \cD\cF_S \to \cD\cAb$ be a triangulated functor and let $M$ be a 1--motive over $S$. There are canonical short exact sequences of $\IZ_\ell$--modules
$$0 \to H^{i-1} F(M)\hotimes \IZ_\ell \to H^i F(\Tell M) \to \Tell H^iF(M) \to 0 $$
and short exact sequences of $\ell$--torsion groups
$$0 \to H^{i-1} F(M)\otimes \IQ_\ell/\IZ_\ell \to H^iF(\Iell M) \to H^iF(M)[\ell^\infty] \to 0 $$
both natural in $M$ and $F$.
\end{prop}

\begin{proof}
The short exact sequence of constant sheaves $0 \to \IZ \xrightarrow{\:\:\ell^i\:\:} \IZ \to \IZ/\ell^i\IZ \to 0$ induces a long exact sequence of groups, from where we can cut out the short exact sequences
$$0 \to H^iF(M) \otimes \IZ/\ell^i\IZ \to H^iF(M \otimes^\IL\IZ/\ell^i\IZ) \to H^{i+1}F(M)[\ell^i] \to 0$$
The limit system of commutative groups $(H^iF(M) \otimes \IZ/\ell^i\IZ)_{i= 0}^\infty$ has the Mittag--Leffler property, and the short exact sequences in the proposition are then obtained by taking limits, respectively colimits over $i \geq 0$.
\end{proof}

\vspace{4mm}
\begin{cor}\label{Cor:TellFixGlobal}
Let $k$ be a number field, let $\ell$ be a prime number and let $M=[u:Y \to G]$ be a 1--motive over $k$. Set $Z := H^{-1}(M) = \ker u$. The morphism of $\IZ_\ell$--modules
$$H^i(k,Z\otimes \IZ_\ell) \to H^i(k,\Tell M)$$
induced by the morphism of 1--motives $[Z\to 0]\to[Y\to G]$ is an isomorphism for $i=0$ and injective for $i=1$.
\end{cor}

\begin{proof}
Proposition \ref{Pro:FiltrationTate} applied to the functor $\IR\Gamma(k,-)$ yields a short exact sequence of $\IZ_\ell$--modules
$$0 \to H^{-1}(k,M) \hotimes \IZ_\ell \to H^0(k,\Tell M) \to \Tell H^0(k,M) \to 0$$
Since $Z(k) = H^{-1}(k,M)$ is a finitely generated group we can identify $Z(k) \hotimes \IZ_\ell$ with $Z(k) \otimes_\IZ \IZ_\ell$, so to get the statement for $i=0$ it remains to show that the last group in this sequence vanishes. Write $\cO_k$ for the ring of integers of $k$, and choose a sufficiently small open subscheme $U\subseteq \spec\cO_k$ such that $M$ extends to a 1--motive over $U$, and such that $\ell$ is invertible on $U$. We have then $Z(U)= Z(k)$ and $H^0(U,\Tell M) = H^0(k, \Tell M)$, so we may as well show that $\Tell H^0(U,M)$ vanishes. Indeed, it follows by d\'evissage from the Mordell--Weil Theorem, Dirichlet's Unit Theorem and finiteness of $H^1(U,Y)$ that $H^0(U,M)$ is a finitely generated group (cf. \cite{Hara05}, Lemma 3.2), so its Tate module is trivial. For the case $i=1$, we consider the triangle
$$[Z\to 0]\to[Y\to G]\to[Y/Z\to G]$$
and observe that if we quotient both terms of the complex $[Y/Z \to G]$ by the finite torsion part of $Y/Z$, we get a quasi--isomorphic complex, which is the complex of a 1--motive $M' = [u':Y' \to G']$ where now $u'$ is injective. By the first part we have $H^0(k,\Tell M')=0$, and the statement can be read in the long exact cohomology sequence associated with $0\to Z\otimes\IZ_\ell \to \Tell M \to \Tell M' \to 0$.
\end{proof}

\vspace{4mm}
\begin{par}
The statement of Corollary \ref{Cor:TellFixGlobal} remains true over any field $k$ which is finitely generated over its prime 
field, and prime number $\ell$ different from the characteristic of $k$. It is wrong in general for local fields. 
\end{par}
\begin{par}
We now come to the dual 1--motive: With each $1$--motive $M$ over a noetherian regular scheme $S$ is functorially associated a dual 1--motive $M^\vee$ over $S$, so that we get an involution of the category $\cM_{1,S}$ of 1--motives over $S$. The duals of tori, lattices and abelian schemes, if seen as a $1$--motives, are the usual duals, and the duality functor is compatible with the weight filtration. This is the content of the following theorem.
\end{par}

\begin{thm}\label{Thm:CartierHauptsatz}
There exists an antiequivalence of categories $(-)^\vee:\cM_{1,S}\to \cM_{1,S}$ such that for every 1--motive $M$ over $S$ the following holds:
\begin{enumerate}
\item[(1)] There are canonical and natural isomorphisms of 1--motives $(M/W_{-i}(M))^\vee \:\:\cong \:\: W_{i-3}(M^\vee)$ for each $i$.
\item[(2)] There is a natural isomorphism $[M^\vee] \cong \IR\cHom(M,\IG_m[1])_{\leq 0}$ in the derived category $\cD\cF_S$, where $(-)_{\leq 0}$ means truncation in degree $0$.
\item[(3)] There is a natural isomorphism of 1--motives $\epsilon_M:M \to M^{\vee\vee}$, such that the morphism of complexes $[\epsilon_M]$ coincides in the derived category of $\cF_S$ with the morphism induced by the canonical evaluation morphism (explained below).
\end{enumerate}
Moreover, the properties \textup{(1)}, \textup{(2)} and \textup{(3)} characterise $(-)^\vee$ up to an isomorphism of functors.
\end{thm}

\begin{par}
For every object $X$ of $\cD\cF_S$ we have a natural morphism $X\to \IR\cHom(\IR\cHom(X,\IG_m[1]),\IG_m[1])$ (see \cite{SGA5}, expos\'e 1, after Proposition 1.6) as well as $X \to X_{\leq 0}$. Together, these yield the natural morphism
$$X_{\leq 0}\to \IR\cHom(\IR\cHom(X,\IG_m[1])_{\leq 0},\IG_m[1])_{\leq 0}$$
which we consider for $X = X_{\leq 0} = [M]$ in part (3) of the theorem.
\end{par}

\vspace{4mm}
\begin{para}\label{Par:CommentsOnDuality}
\begin{par}
The unicity of the functor $(-)^\vee$ can be shown by a simple d\'evissage argument. Its existence is in essence Deligne's construction of the dual 1--motive as given in  \S 10.2.11 of \cite{DeligneHodgeIII}, combined with the following observations (1) and (2): 
\end{par}
\begin{enumerate}
 \item If $X$ is either a finite flat group scheme, a torus or a lattice over $S$, then the sheaf $\cHom(X, \IG_m)$ is represented by the Cartier dual of $X$, and $\cExt^1(X, \IG_m) = 0$. 
 \item If $A$ is an abelian scheme over $S$, the sheaves $\cHom(A,\IG_m)$ and $\cExt^2(A, \IG_m)$ are trivial, and
$\cExt^1(A, \IG_m)$ is represented by the dual abelian scheme $A^\vee$. 
 \item For all $i\geq 2$ the sheaves $\cExt^i(X, \IG_m)$ and $\cExt^i(A, \IG_m)$ are torsion. If $\ell$ is invertible on $S$, these sheaves contain no $\ell$--torsion.
\end{enumerate}
\begin{par}
In the case $X$ is a finite flat group scheme, the statements of (1) can be found in \cite{Oort}, Theorem III.16.1. For locally constant group schemes and tori, these follow from \cite{SGA3} exp. XIII cor. 1.4 and \cite{SGA7} exp. VIII prop. 3.3.1 respectively. In (2), we have $\cHom(A,\IG_m) = 0$ because $A$ is proper and geometrically connected, and $\IG_m$ is affine. The isomorphism $\cExt^1(A, \IG_m) \cong A^\vee$ is given by the classical Barsotti--Weil formula (\cite{Oort}, Theorem III.18.1)\footnote{The additional hypothesis that either $A$ is projective over $S$, or that $S$ is artinian is superfluous. The trouble is caused only by Prop. I.5.3 in \emph{loc.cit.}, where representability of the Picard functor $T \mapsto \Pic A/T$ is known just in these cases. This problem has been overcome by M. Raynaud (\cite{FaltChai90}, Theorem 1.9).\\
Oort proves that the Barsotti--Weil formula over a general scheme holds if it holds over all of its residue fields. He then says that the formula is known to hold over any field and quotes Serre's \emph{Groupes alg\'ebriques et corps de classes}, VII.16, th\'eor\`eme 6. But Serre makes right at the beginning of chapter VII the hypothesis that the ground field is algebraically closed. Hence, as long as all residue fields of $S$ are perfect Oort's proof is fine. The general case follows by checking that Serre's arguments also work verbatim over separably closed fields.}. It is shown in \cite{Bree69} that (over a noetherian regular base scheme, as we suppose $S$ to be) the sheaves $\cExt^i(A, \IG_m)$ are torsion for all $i > 1$. Using the second statement of (1), we see that for $n\neq 0$, the multiplication--by--$n$ map on $\cExt^2(A, \IG_m)$ is injective, hence $\cExt^2(A, \IG_m) = 0$. Finally, the torsion sheaves $\cExt^i(X, \IG_m)$ and $\cExt^i(A, \IG_m)$ contain no $\ell$--torsion, because if $F$ is a finite flat group scheme over $S$ annihilated by $\ell$, then $\cExt^i(F,\IG_m) = 0$ for all $i\geq 1$. Indeed, such a group scheme is locally constant, locally presented as $0 \to \IZ^r \to \IZ^r \to F \to 0$, and the functor $\cHom(\IZ^r,-)$ is exact.
\end{par}
\begin{par}
The reason why we need the truncation operations in part (1) of Theorem \ref{Thm:CartierHauptsatz} is that in general we do not know whether the sheaves $\cExt^i(F, \IG_m)$ vanish for $i >1$ if $F$ is a finite flat group scheme over $S$ which is not locally constant. This is presumably not the case, as an explicit example of Breen suggests (\cite{Bree69b}, Breen works with sheaves for the \'etale topology but it seems that his example also works in the fppf setting). Over a field of characteristic zero, or after inverting all residual characteristics of $S$, the truncation is not needed. 
\end{par}
\end{para}

\vspace{4mm}
\begin{prop}\label{Pro:CartierDualModEll}
Let $M$ be a 1--motive over $S$ with dual $M^\vee$, and let $n\geq 1$ be an integer which is invertible on $S$. The Cartier dual of the finite flat group scheme $H^1([M] \otimes^\IL \IZ/n\IZ)$ is naturally isomorphic to $H^1([M^\vee] \otimes^\IL \IZ/n\IZ)$. In particular, there is a canonical, natural isomorphism
$$\Tell(M^\vee) \cong \cHom(\Tell M, \IZ_\ell(1))$$
of $\ell$--adic sheaves on $S$, for every prime number $\ell$ invertible on $S$.
\end{prop}

\begin{proof}
This follows from Theorem \ref{Thm:CartierHauptsatz} and the statement (3) of \ref{Par:CommentsOnDuality}.
\end{proof}

\vspace{14mm}
\section{\texorpdfstring{The pairing between $\Scha^0(k,M)$ and $\Scha^2(k,M^\vee)$}{The pairing between Scha0(k,M) and Scha2(k,Mv)}}\label{Sec:PairingScha0Scha2}

\begin{par}
We fix number field $k$ with algebraic closure $\overline k$ and write $\Omega$ for the set of all places of $k$. For $v\in \Omega$, we denote by $k_v$ the completion of $k$ at $v$. After recalling the definition of Tate--Shafarevich groups, we use the Poitou--Tate Duality Theorem for finite Galois modules to identify the $\ell$--torsion part of $\Scha^2(k,M^\vee)$ with the dual of $\Scha^1(k, \Tell M)$ for any 1--motive $M=[Y\to G]$ over $k$. Then we show that the group $\Scha^0(k, M)$ is finite, and that its $\ell$--part canonically injects into $\Scha^1(k, \Tell M)$.
\end{par}

\vspace{4mm}
\begin{para}
Let $C$ be a bounded complex of continuous $\Gal(\overline k|k)$--modules. The Tate--Shafarevich groups $\Scha^i(k,C)$ of $C$ are defined by
$$\Scha^i(k, C) := \ker\bigg(H^i(k, C) \to \prod_{v\in \Omega}H^i(k_v, C)\bigg)$$
where $H^i$ is continuous cochain cohomology, with the convention that for archimedean $v$ the group $H^i(k_v, C) = H^i(\Gal(\IC|k_v), C)$ stands for Tate modified cohomology (\cite{Neukirch00}, I\S2). The Tate--Shafarevich groups $\Scha^i(k,M)$ of a 1--motive $M=[Y\to G]$ over $k$ are those of the complex of discrete Galois modules $Y(\overline k) \to G(\overline k)$ placed in degrees $-1$ and $0$.
\end{para}

\vspace{4mm}
\begin{prop}\label{Pro:Scha1TateModScha2}
Let $M$ be a 1--motive over $k$ and let $\ell$ be a prime number. There is a canonical, perfect pairing of topological groups
$$\Scha^1(k, \Tell M)   \:\:\times\:\: \Scha^2(k,M^\vee)[\ell^\infty] \:\:\to \QZ$$
where $\Scha^1(k, \Tell M)$ is profinite and $\Scha^2(k,M^\vee)[\ell^\infty]$ is discrete. In particular, $\Scha^2(k,M^\vee)[\ell^\infty]$ is finite or zero if and only if $\Scha^1(k, \Tell M)$ is so.
\end{prop}

\begin{proof}
By Poitou--Tate duality for finite Galois modules (\cite{Neukirch00}, Theorem 8.6.8) we have a natural, perfect duality between finite groups 
$$\Scha^1(k,M\otimes^\IL \IZ/\ell^i\IZ)\times \Scha^2(k,M^\vee\otimes^\IL \IZ/\ell^i\IZ) \to \QZ$$
noting Proposition \ref{Pro:CartierDualModEll}. The functor $\Scha^1(k, -)$ commutes with limits of finite Galois modules by \cite{Serre1} Proposition 7, and that $\Scha^2(k,-)$ commutes with arbitrary colimits, so we obtain a perfect pairing of topological groups 
$$\Scha^1(k, \Tell M) \times \Scha^2(k,\Iell M^\vee) \to \IQ/\IZ$$
It remains to show that $\Scha^2(k, \Iell M^\vee)$ is canonically isomorphic to the $\ell$--part of the torsion group $\Scha^2(k, M^\vee)$. Indeed, from Proposition \ref{Pro:FiltrationTate} we get the following commutative diagram of torsion groups with exact rows
$$\begin{diagram}
\setlength{\dgARROWLENGTH}{4mm}%
\node{0}\arrow{e}\node{H^1(k,M^\vee)\otimes \IQ_\ell/\IZ_\ell}\arrow{s}\arrow{e}\node{H^2(k,\Iell M^\vee)}\arrow{s} \arrow{e} \node{H^2(k, M^\vee)[\ell^\infty]}\arrow{s}\arrow{e}\node{0}\\
\node{0}\arrow{e}\rnode{\!\!\prod_{v\in \Omega}\!H^1(k_v, M^\vee)\!\otimes\!\IQ_\ell/\IZ_\ell\!\!}\arrow{e}\rnode{\!\!\prod_{v\in \Omega}\! H^2(k_v, \Iell M^\vee)\!\!}\arrow{e}\rnode{\!\prod_{v\in \Omega}\! H^2(k_v, M^\vee)[\ell^\infty]\!\!}\arrow{e}\node{0}
\end{diagram}$$
Because $H^1(k,M^\vee)$ and $H^1(k_v,M^\vee)$ are torsion, the first terms of both rows are zero, hence the canonical isomorphism $\Scha^2(k,\Iell M^\vee) \cong \Scha^2(k,M^\vee)[\ell^\infty]$ as required.
\end{proof}

\vspace{4mm}
\begin{para}\label{Par:CanonicalAst}
Let $M$ be a 1--motive over $k$ and let $\ell$ be a prime number. From Proposition \ref{Pro:FiltrationTate} we get a commutative diagram of $\IZ_\ell$--modules with exact rows
$$\begin{diagram}
\setlength{\dgARROWLENGTH}{4mm}%
\node{0}\arrow{e}\node{H^0(k,M)\hotimes\IZ_\ell}\arrow{s,l}{\alpha_\ell}\arrow{e}\node{H^1(k,\Tell M)} \arrow{s} \arrow{e} \node{\Tell H^1(k, M)}\arrow{s}\arrow{e}\node{0}\\
\node{0}\arrow{e}\rnode{\!\!\prod_{v\in \Omega}\!H^0(k_v,M)\!\hotimes\!\IZ_\ell\!\!}\arrow{e}\rnode{\!\!\prod_{v\in \Omega} H^1(k_v, \Tell M)\!\!}\arrow{e}\rnode{\!\prod_{v\in \Omega} \Tell H^1(k_v, M)\!\!}\arrow{e}\node{0}
\end{diagram}$$
The kernel of the right hand vertical map is the Tate module of $\Scha^1(k,M)$, which is torsion free, and even trivial if $\Scha^1(k,M)$ is finite (which conjecturally is always so, c.f.\ \cite{Hara05}, Corollary 4.9). In any case the map $\ker\alpha_\ell \to \Scha^1(k,\Tell M)$ is an isomorphism on torsion elements. The kernel of $\alpha_\ell$ contains $\Scha^0(k,M)\otimes \IZ_\ell$, hence an injection
\begin{equation}
\Scha^0(k,M)\otimes \IZ_\ell \to \Scha^1(k,\Tell M)\tag{$\ast$}
\end{equation}
In Section 5 of \cite{Hara05} a profinite group $\Scha^0_\wedge(k,M)$ was introduced. Its pro--$\ell$ part is $\ker\alpha_\ell$, hence equal to $\Scha^1(k,\Tell M)$ in the case $\Scha^1(k,A)$ is finite. This relates Proposition 5.1 of \emph{loc.cit.}\ to our Proposition \ref{Pro:Scha1TateModScha2}. There is a canonical isomorphism $\Scha^0(k,M)\otimes \IZ_\ell \cong \Scha^0(k,M)[\ell^\infty]$ because $\Scha^0(k,M)$ is finite as we shall see in Proposition \ref{Pro:CalculScha0}. These observations yield the following corollary to Proposition \ref{Pro:Scha1TateModScha2}.
\end{para} 

\vspace{4mm}
\begin{cor}\label{Cor:Scha0Scha2PairingDef}
The pairing of proposition \ref{Pro:Scha1TateModScha2} induces a pairing
$$\Scha^0(k,M)[\ell^\infty] \times \Scha^2(k,M^\vee)[\ell^\infty] \to \QZ$$
which is nondegenerate on the left. Its right kernel is divisible if and only if the map $(\ast)$ induces an isomorphism $\Scha^0(k,M)\otimes \IZ_\ell \to \Scha^1(k,\Tell M)_\uptor$, and this pairing is a perfect pairing of finite groups if and only if the map $(\ast)$ is an isomorphism. 
\end{cor}

\vspace{4mm}
\begin{par}
We end the section with the following proposition which explains the group $\Scha^0(k,M)$ and shows that it is finite (cf. \cite{Hara05} Lemma 4.11). 
\end{par}

\begin{prop}\label{Pro:CalculScha0}
Let $M = [u:Y\to G]$ be a 1--motive over $k$ and set $Z := \ker u$. The morphism of 1--motives $[Z\to 0]\to[Y\to G]$ induces an isomorphism of finite groups
$$\Scha^1(k,Z) \xrightarrow{\:\:\cong\:\:} \Scha^0(k,M)$$
Moreover, for any prime number $\ell$, there are canonical isomorphisms of finite groups
$$\Scha^1(k,Z)[\ell^\infty] \cong \Scha^1(k,Z)\otimes \IZ_\ell \cong \Scha^1(k,Z\otimes \IZ_\ell)\cong \Scha^1(k, \Tell[Z \to 0])$$
All these groups are annihilated by the order of any finite Galois extension $k'|k$ over which $Z$ is constant.
\end{prop}

\begin{proof}
By diagram chase, using one finite place $v\in \Omega$, we see that the map $\Scha^0(k,M) \to \Scha^1(k,Y)$ is injective. It follows in particular that $\Scha^0(k,M)$ is zero if the Galois action on $Y$ is trivial. In general, let $k'|k$ be a finite Galois extension such that $\Gal(\overline k|k')$ acts trivially on $Y$, and let $\Omega'$ be the set of places of $k'$. For $w \in \Omega'$ lying over $v\in \Omega$ we write $k_w'$ for the completion of $k'$ at $w$ and $k_w$ for the completion of $k$ at $v$. As a Galois module, $Z := \ker(u)$ can be interpreted as $Z = H^{-1}(k',M)$. From the Hochschild--Serre spectral sequence we get then a commutative diagram with exact rows
$$\begin{diagram}
\setlength{\dgARROWLENGTH}{4mm}
\node{0}\arrow{e}\node{H^1(\Gal(k'|k),Z)} \arrow{s}\arrow{e}\node{H^0(k,M)}
\arrow{s}\arrow{e} \node{H^0(k',M)} \arrow{s}\\
\node{0}\arrow{e} \rnode{\prod_{w\in \Omega'} H^1(\Gal(k'_w|k_w),Z)} \arrow{e} \rnode{\prod_{w\in \Omega'}
H^0(k_w,M)} \arrow{e} \rnode{\prod_{w\in \Omega'} H^0(k'_w,M)}
\end{diagram}$$
Because $\Gal(\overline k|k')$ acts trivially on $Y$ we have $\Scha^0(k',M)=0$ by our previous observation, hence
$$\Scha^0(k,M) \cong \ker\left(H^1(\Gal(k'|k),Z) \to \prod_{w\in \Omega'}H^1(\Gal(k'_w|k_w),Z)\right)$$
the product running over all $w\in \Omega'$, or alternatively, over all decomposition subgroups of $\Gal(k'|k)$. The finiteness statement follows, as $H^1(\Gal(k'|k),Z)$ is finite and annihilated by the order of $\Gal(k'|k)$ (\cite{Weib94}, Theorem 6.5.8 and Corollary 6.5.10). Repeating the arguments for the 1--motive $[Z\to 0]$ yields the first isomorphism of the proposition.\\
Now let $\ell$ be a prime number. The first isomorphism from the left exists for any finite commutative group in place of $\Scha^1(k,Z)$. For the next one, choose a finite Galois extension $k'|k$ such that $Z$ is constant over $k'$. We can proceed as before in order to express $\Scha^1(k,Z)$ and $\Scha^1(k,Z\otimes \IZ_\ell)$ in terms of cohomology groups of the finite group $\Gal(k'|k)$ and its subgroups. It remains to show that, given a finite group $\Gamma$ acting on $Z$, the canonical map $H^1(\Gamma,Z)\otimes\IZ_\ell \to H^1(\Gamma,Z \otimes\IZ_\ell)$ is an isomorphism. This is indeed so for any flat $\IZ$--algebra in place of $\IZ_\ell$ by the universal coefficient theorem. The last isomorphism holds because $Z\otimes\IZ_\ell \cong \lim Z/\ell^iZ$ and because limits are left exact and commute with continuous $H^1$. 
\end{proof}

\vspace{14mm}
\section{Lie algebra cohomology of the Tate module}\label{Sec:H1ast}

\begin{par}
We fix a number field $k$ with algebraic closure $\overline k$ and a prime number $\ell$. With every 1--motive $M$ over $k$ is associated a continuous $\IQ_\ell$--linear representation $V_\ell M = \Tell M \otimes \IQ_\ell$ of $\Gal(\overline k|k)$. The image of $\Gal(\overline k|k)$ in $\GL(\Vell M)$ is an $\ell$--adic Lie group $L^M$, whose Lie algebra we denote by $\fl^M$. An idea going back to Serre and Tate, used by Serre in \cite{Serre1} to solve the congruence subgroup problem for abelian varieties over number fields, is to consider the vector space $H^1_\ast(\fl^M,\Vell M)$ consisting of those elements of $H^1(\fl^M,\Vell M)$ which restrict to zero on each one dimensional subalgebra of $\fl^M$. Our goal is to describe $H^1_\ast(\fl^M,\Vell M)$. 
\end{par}

\begin{par}
We will work only with 1--motives $M=[Y\to G]$ where $G$ is either an abelian variety or a torus, rather than a general semiabelian variety. This brings considerable simplifications in both, statements and proofs. I will comment at the end of the section (\ref{Para:CommentOnGeneralG}) on this hypothesis and on the modifications that are necessary in order to compute $H^1_\ast(\fl^M,\Vell M)$ for general 1--motives. The following theorem is the crucial ingredient for our finiteness results.
\end{par}

\begin{thm}\label{Thm:LieCohoMain}
Let $M=[u := Y\to G]$ be a 1--motive over $k$ where $G$ is an abelian variety or a torus. Set $E_\ell := \End_{\overline k}(G) \otimes\IQ_\ell$ and $X_\ell := \im(u)\otimes\IQ_\ell$, denote by $D_\ell$ the $E_\ell$--submodule of $G(\overline k)\otimes\IQ_\ell$ generated by $X_\ell$, and define 
$$\overline X_\ell := \big\{x\in D_\ell \:\big|\: f(x)\in f(X_\ell) \pt f\in\Hom_{E_\ell}(D_\ell, \Vell G) \big\}$$
There is a canonical isomorphism of $\IQ_\ell$--vector spaces $\overline X_\ell/X_\ell \cong H^1_\ast(\fl^M,\Vell M)$.
\end{thm}

\begin{par}
\noindent The proof of this theorem relies on a structure result for the Lie algebra $\fl^M$ which in turn relies on Faltings's theorems on endomorphisms of abelian varieties over number fields. Observe that the object $\overline X_\ell/X_\ell$ can be calculated by means of ordinary linear algebra. The statement of the theorem is wrong for general semiabelian varieties $G$.
\end{par}

\vspace{4mm}
\begin{para}
\begin{par}
We recall some definitions and results from \cite{Serre1}. Let $L$ be a profinite group and let $T$ be a continuous $L$--module. We write $H^1(L,T)$ for the group of continuous cocycles $L\to T$ modulo coboundaries and define
$$H^1_\ast(L,T) := \ker\bigg(H^1(L,T) \to \prod_{x\in L}H^1(\angl{x},T)\bigg)$$
where $\angl{x}$ denotes the closed subgroup of $L$ generated by $x$. If $N$ is a closed normal subgroup of $L$ acting trivially on $T$, then the inflation map induces an isomorphism $H^1_\ast(L/N,T)\to H^1_\ast(L,T)$ (\emph{loc.cit.}, Proposition 6). If $T$ is a profinite $L$--module, say $T=\lim T_i$ where the $T_i$ are finite discrete $L$--modules, then the canonical map $H^1(L,T) \to \lim H^1(L,T_i)$ is an isomorphism (\emph{loc.cit.}, Proposition 7). Because the limit functor is right exact, also the  canonical map $H^1_\ast(L,T) \to \lim H^1_\ast(L,T_i)$ is an isomorphism in that case.
\end{par}
\begin{par}
For a Lie algebra $\fl$ acting on a vector space $V$, we denote by $H^1_\ast(\fl,V)$ the subspace of $H^1(\fl,V)$ consisting of those elements which restrict to zero in $H^1(\angl x,V)$ for every 1--dimensional subalgebra $\angl x$ of $\fl$.
\end{par}
\end{para}

\vspace{4mm}
\begin{lem}\label{Lem:LazardH1ast}
Let $L$ be a compact $\ell$--adic Lie group with Lie algebra $\fl$, acting on a finite dimensional $\IQ_\ell$--vector space $V$. For any open subgroup $N$ of $L$, the equality 
$$H^1_\ast(L,V) = \ker\bigg(H^1(L,V) \to \prod_{x\in N}H^1(\angl{x}, V) \bigg)$$
holds. Moreover, if $N$ is normal there is a canonical isomorphism $H^1_\ast(L,V) \cong H^1_\ast(N,V)^{L/N}$, and if $N$ is sufficiently small, there is a canonical isomorphism $H^1_\ast(N,V) \cong H^1_\ast(\fl, V)$.
\end{lem}

\begin{proof}
Let $N$ be an open subgroup of $L$, and let $c$ be an element of $H^1(L,V)$ restricting to zero in $H^1(\angl{x},V)$ for each $x\in N$. Fix an element $x\in L$ and let us show that $c$ restricts to zero in $H^1(\angl{x},V)$. Because $\angl{x}$ is compact, the quotient $\angl{x}/(N\cap \angl{x})$ is finite. By a restriction--corestriction argument and using that $V$ is uniquely divisible we see that the restriction map $H^1(\angl{x},V) \to H^1(\angl{x}\cap N,V)$ is injective, hence the first claim. Now suppose that $N$ is open and normal. Since $L$ is compact, the quotient $L/N$ is finite and we have $H^i(L/N,V) = 0$ for all $i>0$, and the Hochschild--Serre spectral sequence yields an isomorphism $H^1(L,V) \cong H^1(N,V)^{L/N}$. We must show that in the diagram
$$\begin{diagram} \setlength{\dgARROWLENGTH}{4mm} 
\node{0}\arrow{e}\node{H^1_\ast(L, V)}\arrow{s}\arrow{e}\node{H^1(L, V)}\arrow{s,l}{\cong}\arrow{e}\arrow{se,b,--}{\delta} \rnode{\prod_{x\in L}H^1(\angl{x}, V)}\arrow{s}\\
\node{0}\arrow{e}\node{H^1_\ast(N, V)^{L/N}}\arrow{e}\node{H^1(N,V)^{L/N}\!\!\!}\arrow{e}\rnode{\!\!\prod_{x\in N}H^1(\angl{x},V)}
\end{diagram}$$
the left hand vertical map is an isomorphism, i.e.\ that the kernel of the diagonal map $\delta$ is exactly $H^1_\ast(L, V)$. But this is again the first statement of the lemma. Finally, if $N$ is sufficiently small we have an isomorphism $H^1(N, V) \cong H^1(\fl, V)$ by a well--known theorem of Lazard (\cite{Lazard}, V.2.4.10), from which the last statement follows.
\end{proof}

\vspace{4mm}
\begin{para}\label{Par:flMIntro}
Let $M=[u:Y\to G]$ be a 1--motive over $k$ where $G$ is an abelian variety or a torus. The Tate module $\Tell M$ of $M$ is an extension of $Y\otimes \IZ_\ell$ by the Tate module $\Tell G$ of $G$ as we have seen in \ref{Par:TellMExplicit}, so we get an extension continuous of Galois representations
$$0 \to \Vell G \to \Vell M \to Y\otimes\IQ_\ell \to 0$$
We denote by $\fl^M$ and $\fl^G$ the Lie algebra of the image of $\Gamma := \Gal(\overline k|k)$ in the group $\GL(\Vell M)$, respectively in $\GL(\Vell G)$. The Galois group $\Gamma$ acts continuously on these Lie algebras by conjugation, and we have a canonical surjection $\fl^M\to \fl^G$. Let $\fn^M$ denote its kernel, so $\fn^M$ consists of those elements of $\fl^M$ which act trivially on $\Vell G$. The Lie algebra $\fn^M$ is commutative, and we can identify it with a linear subspace of $\Hom(Y\otimes\IQ_\ell,\Vell G)$ via the map 
$$\vartheta:\fn^M \to \Hom_{\IQ_\ell}(Y\otimes\IQ_\ell,\Vell G)$$
given by $\vartheta(x)(y) = x.v$ where $v\in\Vell M$ is any element mapping to $y \in Y\otimes\IQ_\ell$. Routine checking shows that this map is well defined, injective and $\Gamma$--equivariant. We can describe the image of $\vartheta$ as follows: Look at $u$ as being a $k$--rational point on the abelian variety or torus $\cHom(Y,G)$, and denote by $B$ the connected component of the smallest algebraic subgroup of $\cHom(Y,G)$ containing $u$. Then $B$ is also an abelian variety or a torus, and we have an inclusion
$$\Vell B \subseteq \Vell\cHom(Y, G)\cong\Hom_{\IQ_\ell}(Y\otimes\IQ_\ell,\Vell G)$$
The following theorem is a special case of Theorem 6.2 of \cite{Jos11}. In the case $G$ is an abelian variety it goes back to a result of Ribet (\cite{Ribet}, see Appendix 2 of \cite{Hind88}).
\end{para}

\begin{thm}\label{Thm:Ribet}
The map $\vartheta$ induces an isomorphism of Galois representations $\vartheta: \fn^M \xrightarrow{\:\:\cong\:\:} \Vell B$.
\end{thm}

\begin{par}
In particular, it follows that the dimension of $\fn^M$ is independent of $\ell$. If $G$ is an abelian variety, it is still unknown whether the dimension of $\fl^G$ is independent of $\ell$.
\end{par}

\vspace{4mm}
\begin{lem}\label{Lem:DandHomBG}
Let $M = [u:Y\to G]$ be a 1--motive over $k$ where $G$ is an abelian variety or a torus. Denote by $D$ the $\End_{\overline k}(G)$ submodule of $G(\overline k)$ generated by $\im(u)$ and define $B\subseteq \cHom(Y,G)$ as in \ref{Par:flMIntro}. The linear map
$$h:\Hom_{\overline k}(B,G) \otimes\IQ \to G(\overline k)\otimes\IQ$$
given by $h(\psi \otimes 1) = \psi(nu) \otimes n^{-1}$, where $n\geq 1$ is any integer such that $nu\in B(k)$, induces an isomorphism $\Hom_{\overline k}(B,G) \otimes\IQ \cong D\otimes\IQ$.
\end{lem}

\begin{proof}
The homomorphism $h$ is injective. Indeed, if $\psi:B\to G$ is such that $h(\psi)=0$, then $\ker\psi$ is a subgroup of $B$ containing a nonzero multiple of $u$, hence $\ker\psi = B$ by minimality of $B$. By Poincar\'e's Complete Reducibility Theorem (\cite{Mumf70} IV.19, Theorem 1), the inclusion $B\subseteq \cHom(Y,G)$ induces a surjection
$$Y \otimes \End_{\overline k}(G) \cong \Hom_{\overline k}(\cHom(Y,G),G) \xrightarrow{\:\:\res\:\:} \Hom_{\overline k}(B,G)$$
sending $y\otimes \varphi$ to the unique homomorphism $\psi:B\to G$ with $\psi(nu) = n\varphi(u(y))$ where $n\geq 1$ is sufficiently big to that $nu\in B(k)$, so the remaining statements follow.
\end{proof}

\vspace{4mm}
\begin{lem}\label{Lem:CommutativeVW}
Let the 1--motive $M = [Y\to G]$, the subgroup $D\subseteq G(\overline k)$ and the algebraic subgroup $B \subseteq G$ be as in Lemma \ref{Lem:DandHomBG}. There is a commutative diagram
$$\begin{diagram}
\setlength{\dgARROWLENGTH}{4mm}
\node{Y \otimes \IQ_\ell}\arrow{s,r}{\cong}\arrow{e,t}{u\otimes \id}\node{D \otimes \IQ_\ell}\arrow{s,lr}{(2)}{\cong} \arrow{e,tb}{(1)}{\cong}\node{\Hom_{\fl^G}(\Vell B, \Vell G)}\arrow{s,lr}{\cong}{\vartheta^\ast}\\
\node{H^0(\fl^M, Y\otimes \IQ_\ell)}\arrow{e,t}{\partial}\node{H^1(\fl^M, \Vell G)}\arrow{e,tb}{\res}{\cong}\node{\Hom_{\fl^G}(\fn^M, \Vell G)}
\end{diagram}$$
with canonical isomorphisms where indicated.
\end{lem}

\begin{proof}
\begin{par}
We start wit the left hand square. The left hand vertical isomorphism is tautological, because $\fl^M$ acts trivially on $Y\otimes \IQ_\ell$. The map $\partial$ is the connecting morphism in the long exact cohomology sequence coming from the weight filtration of $M$. The vector spaces $D\otimes \IQ_\ell$ and $H^1(\fl^M, \Vell G)$ are naturally $E_\ell := \End_{\overline k}(G)\otimes\IQ_\ell$--modules, the first one by definition, and the second one via the canonical action of $E_\ell$ on $\Vell G$. The map (2) is then given by $E_\ell$--linearity and sending $u(y)\otimes 1$ to $\partial(y \otimes 1)$ for all $y\in Y$. By definition of $D$ this indeed describes a unique map such that the left hand square commutes. We will see in a moment that it is well defined and an isomorphism.
\end{par}
\begin{par}
We now come to the right hand square, starting with the description of the map (1). Every element of $D\otimes\IQ_\ell$ is a linear combination of elements of the form $\psi(u)\otimes 1$ for some $\psi\in\Hom_{\overline k}(B,G)$ by Lemma \ref{Lem:DandHomBG}. The map (1) is given by linearity and sends $\psi(u) \otimes 1$ to the $\fl^G$--equivariant map $\Vell \psi: \Vell B \to \Vell G$. This map is an isomorphism by Lemma \ref{Lem:DandHomBG} and by Faltings's theorem on homomorphisms of abelian varieties. The right hand vertical map is given by precomposition with the isomorphism $\vartheta$, hence an isomorphism. The lower horizontal map is given by restriction of cocycles, and an isomorphism because $H^i(\fl^G,\Vell G)$ vanishes for $i=1,2$ (c.f. \cite{Serre2} Th\'eor\`eme 2). 
\end{par}
\begin{par}
By definition of $\vartheta$, the big square commutes. Moreover, the isomorphisms (1), $\vartheta^\ast$ and the inverse of $\res$ are all isomorphisms of $E_\ell$--modules. Hence so is their composition, which is then an isomorphism of $E_\ell$--modules $D\otimes \IQ_\ell \to H^1(\fl^M, \Vell G)$ which must coincide with (2). 
\end{par}
\end{proof}

\vspace{4mm}
\begin{prop}\label{Pro:flMSections}
The Lie algebra extension $0\to \fn^M \to \fl^M \to \fl^G\to 0$ is split. There exist a Lie algebra section $\sigma:\fl^G \to \fl^M$ and a $\IQ_\ell$--linear section $s:Y\otimes \IQ\to \Vell M$ such that the action of $\fl^M$ on $\Vell M$ is given by
$$(f+\sigma(g)).(v + s(y)) = g.v + f.y$$
for all $f\in \fn^M$, all $g\in \fl^G$, all $v\in \Vell A$ and $y\in Y\otimes \IQ_\ell$.
\end{prop}

\begin{proof}
This is Corollary 2.19 of \cite{Jos10}. It is essentially a consequence of Theorem \ref{Thm:Ribet}, semi\-simplicity of $\Vell G$ as $\fl^G$--module (Faltings, \cite{Faltings}) and the vanishing of $H^i(\fl^G, \Vell G)$ for $i= 1,2$ (\cite{Serre2} Th\'eor\`eme 2). 
\end{proof}

\vspace{4mm}
\begin{lem}\label{Lem:flMSections}
Let $M=[u := Y\to G]$ be a 1--motive over $k$ where $G$ is an abelian variety or a torus. In order that an element $h \in H^1(\fl^M, \Vell G)$ belongs to $H^1_\ast(\fl^M, \Vell M)$ it suffices that it maps to zero in $H^1(\fc, \Vell M)$ for each one dimensional subalgebra $\fc$ of $\fn^M$.
\end{lem}

\begin{proof}
Represent $h \in H^1(\fl^M, \Vell G)$ by a cocycle $c:\fl^M\to \Vell G \subseteq \Vell M$, and choose a linear section $s:Y\otimes\IQ_\ell \to \Vell M$ and a Lie algebra section $\sigma:\fl^G\to\fl^M$ as in \ref{Pro:flMSections}. Since $H^1(\fl^G,\Vell G)$ vanishes (\cite{Serre2} Th\'eor\`eme 2) the cocycle $c\circ\sigma$ is a coboundary. Thus, changing $c$ by a coboundary we may suppose that $c\circ\sigma
=0$. Let $\fc$ be a 1--dimensional subalgebra of $\fl^M$ generated by an element $x\in\fl^M$. We have to show that there exists an element $v\in \Vell M$ such that $c(x)=x.v$. We can write $x$ as $x=f+\sigma(g)$ for some $f\in\fn^M$ and $g\in\fl^G$. By hypothesis, there exists an element $v\in\Vell M$ with $c(f)=f.v$. We can write $v$ as $v=v'+s(y)$ for some $v'\in\Vell G$ and $y\in Y\otimes \IQ_\ell$. We then have 
$$c(x) = c(f+\sigma(g)) = c(f) = f.v = f.s(y) = (\sigma(g)+f).s(y) = x.s(y)$$
and this proves the lemma.
\end{proof}

\vspace{4mm}
\begin{lem}\label{Lem:KernelCorresp}
Let $M = [u:Y\to G]$ be a 1--motive over $k$ where $G$ is an abelian variety or a torus. Denote by $D$ the $E:=\End_{\overline k}(G)$ submodule of $G(\overline k)$ generated by $X := \im(u)$. The isomorphism $D\otimes\IQ_\ell\to H^1(\fl^M, \Vell G)$
from Lemma \ref{Lem:CommutativeVW} induces an isomorphism between the kernels of the maps
$$D\otimes \IQ_\ell \to \prod_{h\in H} \Vell G/h(X\otimes\IQ_\ell) \qqet H^1(\fl^M, \Vell G)\to \prod_{\fc \subseteq \fl^M} H^1(\fc, \Vell M)$$
where the left hand product runs over all $h \in H := \Hom_{E\otimes\IQ_\ell}(D \otimes\IQ_\ell,\Vell G)$.
\end{lem}

\begin{proof}
Lemma \ref{Lem:flMSections} shows that if on the right hand side we let the product only run over $\fc \in \fn^M$ we still get the same kernel. For every $\fc = \angl{x} \subseteq \fn^M$ we have
$$H^1(\fc, \Vell M) \:\cong\: \frac{\Vell M}{\{x.v\tq v\in \Vell M\}} \:=\: \frac{\Vell M}{\im(\vartheta(x))}$$
where $\vartheta:\fn\to \Hom(Y\otimes\IQ,\Vell G)$ is as defined in \ref{Par:flMIntro}. The map $H^1(\fl^M, \Vell G) \to \Hom_{\fl^G}(\fn^M, \Vell G)$ given by restriction of cocycles is an isomorphism, thus we have to show that the kernels of the maps
$$D \otimes \IQ_\ell \to \prod_{h\in H} \Vell G/h(X\otimes\IQ_\ell) \qqet \Hom_{\fl^G}(\fn^M, \Vell G) \to \prod_{x \in \fn^M} \Vell G/\im(\vartheta(x))$$
correspond under the isomorphism $D \otimes \IQ_\ell \cong \Hom_{\fl^G}(\fn^M, \Vell G)$ sending $\psi(u)\otimes 1$ to $\Vell\psi\circ\vartheta$ for all $\psi\in\Hom_{\overline k}(B,G)$. Here, $B\subseteq \cHom(Y,G)$ is defined as in Theorem \ref{Thm:Ribet}. The right hand map  sends an $\fl^G$---module homomorphism $c: \fn^M \to \Vell G$ to the class of $c(x)$ in the factor corresponding to $x$. By Lemma  \ref{Lem:CommutativeVW}, we can on the left hand side as well take $\Hom_{\fl^G}(\Vell B, \Vell G)$ in place of $D \otimes \IQ_\ell$. Then we must show that the kernels of the maps
$$\Hom_{\fl^G}(\Vell B, \Vell G) \to \prod_h \Vell G/f(X\otimes \IQ_\ell) \qqet \Hom_{\fl^G}(\fn^M, \Vell G) \to \!\!\! \prod_{x\in \fn^M }\!\!\! \Vell G/\im(\vartheta(x))$$
correspond to each other via composition with the isomorphism $\vartheta: \fn^M \to \Vell B$, the left hand product now running over all $E\otimes\IQ_\ell$--module morphisms $h:\Hom_{\fl^G}(\Vell B, \Vell G)\to \Vell G$. The canonical map
$$\Vell B \xrightarrow{\:\:\cong\:\:} \Hom_{E \otimes \IQ_\ell}(\Hom_{\fl^G}(\Vell B, \Vell G),\Vell G) \qquad v\mapsto[f\mapsto f(v)]$$
is an isomorphism by Schur's Lemma, so all these $E \otimes \IQ_\ell$--module homomorphisms $h$ are given by evaluation in an element $v\in \Vell B$. If $h$ is the evaluation in $v = \vartheta(x)$ for some $x \in \fn^M$, then $h(X\otimes \IQ_\ell) = \im\vartheta(x)$, hence the claim of the lemma.
\end{proof}

\vspace{4mm}
\begin{proof}[Proof of Theorem  \ref{Thm:LieCohoMain}]
We consider the following diagram, where the exact row is induced by the weight filtration on the $\fl^M$--module $\Vell M$, and where the column is exact by definition:
$$\begin{diagram}
\setlength{\dgARROWLENGTH}{4mm}
\node[3]{H^1_\ast(\fl^M,\Vell M)}\arrow{s,r}{\subseteq}\arrow{se,t}{0}\\
\node{H^0(\fl^M,Y\otimes \IQ_\ell)}\arrow{e,t}{\partial}\node{H^1(\fl^M,\Vell G)}\arrow{se,b}{\delta}\arrow{e}\node{H^1(\fl^M,\Vell M)} \arrow{s}\arrow{e} \node{H^1(\fl^M, Y\otimes\IQ_\ell)}\\
\node[3]{\prod_{\fc \subseteq \fl^M}H^1(\angl{x},\Vell M)}
\end{diagram}$$
The upper diagonal map is zero, because $\fl^M$ acts trivially on $Y\otimes \IQ_\ell$, hence $H^1_\ast(\fl^M, Y\otimes\IQ_\ell)$ is trivial. This shows that every element of $H^1_\ast(\fl^M,\Vell M)$ comes from an element in $H^1(\fl^M,\Vell G)$, so we find an isomorphism 
$$\ker\delta/\im\partial \cong H^1_\ast(\fl^M, \Vell M) $$
induced by the inclusion $\Vell G \subseteq \Vell M$. The lemmas \ref{Lem:CommutativeVW} and \ref{Lem:KernelCorresp} respectively show that the isomorphism $D \otimes \IQ_\ell\to H^1(\fl^M, \Vell G)$ induces isomorphisms $X \otimes\IQ_\ell \cong \im\partial$ and 
$$\overline X_\ell := \{x\in D \otimes\IQ_\ell\tq f(x) \in f(X\otimes\IQ_\ell)\pt f\in \Hom_{E\otimes\IQ_\ell}(D \otimes \IQ_\ell, \Vell G)\} \cong \ker\delta$$
as needed.
\end{proof}

\vspace{4mm}
\begin{para}\label{Para:CommentOnGeneralG}
\begin{par}
Throughout this section we have always supposed that the semiabelian variety $G$ is either an abelian variety or a torus. Most statements and constructions, notably Theorems \ref{Thm:LieCohoMain} and \ref{Thm:Ribet}, remain true if $G$ is isogenous to a product of an abelian variety and a torus, and the proofs require only small additional arguments, but the statements are wrong for general semiabelian varieties. The main problem here is that a general semiabelian variety $G$ is not a semisimple object, so the analogue of Poincar\'e's Complete Reducibility Theorem fails, and the Galois representation $\Vell G$ is not semisimple either.
\end{par}
\begin{par}
In a general setting, the Lie algebra $\fn^M$ should be replaced by the subalgebra of $\fl^M$ consisting of those elements of which act trivially on $\Vell A$ and $\Vell T$, where $A$ and $T$ are the abelian, respectively torus part of $M$. This is then in general not a commutative, but just a nilpotent Lie algebra. The generalisation of Theorem \ref{Thm:Ribet} is Theorem 6.2 in \cite{Jos11}. The subgroup $D$ of $G(\overline k)$ has to be replaced by the group of so--called \emph{deficient points} (\emph{loc.cit.}\ Definition 6.2), and the generalisation of Lemma \ref{Lem:DandHomBG} is \emph{loc.cit.} Theorem 8.10. Finally, $E$--linearity should be reformulated in terms of derivations. With these settings, it should be possible to generalise Theorem \ref{Thm:LieCohoMain} to general 1--motives.
\end{par}
\end{para}

\vspace{14mm}
\section{Finiteness results}\label{Sec:Finiteness}

\begin{par}
In this section we prove the finiteness statements of Theorem \ref{Thm:Main3} stated in the introduction. We fix a number field $k$ with algebraic closure $\overline k$, a prime number $\ell$ and write $\Gamma := \Gal(\overline k|k)$ for the absolute Galois group of $k$ and $\Omega$ for the set of all places of $k$. For a 1--motive $M$ over $k$, we write $\Vell M := \Tell M \otimes\IQ_\ell$, denote by $L^M$ the image of $\Gamma$ in $\GL(\Tell M)$ and by $\fl^M \subseteq \End(\Vell M)$ the Lie algebra of $L^M$.
\end{par}

\vspace{4mm}
\begin{thm}\label{Thm:FinitenessEssentials}
Let $M=[u:Y\to G]$ be a 1--motive over $k$. The $\IZ_\ell$--module $\Scha^1(k,\Tell M)$ is finitely generated. If the semiabelian variety $G$ is a an abelian variety or a torus such that $\End_{\overline k}(G)\otimes\IQ$ is a product of division algebras, then $\Scha^1(k,\Tell M)$ is finite.
\end{thm}

\vspace{4mm}
\begin{par}
The plan of this section is as follows: First we show that for every 1--motive $M$ over $k$ there is a canonical injection of $\Scha^1(k,\Tell M)$ into $H^1_\ast(L^M, \Tell M)$. We continue with some elementary linear algebra and prove, using Theorem \ref{Thm:LieCohoMain}, that the group $H^1_\ast(L^M, \Tell M)$, and hence $\Scha^1(k,\Tell M)$, is finite for all $\ell$ if $M$ is a 1--motive satisfying the condition in the theorem.
\end{par}

\vspace{4mm}
\begin{prop}\label{Pro:SchaAndH1astRank}
Let $M$ be a 1--motive over $k$. There is a canonical injective $\IZ_\ell$--linear map $\Scha^1(k,\Tell M) \to H^1_\ast(L^M, \Tell M)$. The $\IZ_\ell$--module $H^1_\ast(L^M, \Tell M)$ is finitely generated and its rank is bounded by the dimension of $H^1_\ast(\fl^M, \Vell M)$.
\end{prop}

\begin{proof}
For every finite Galois module $F$ the subgroup $\Scha^1(k,F)$ of $H^1(k,F) = H^1(\Gamma,F)$ is contained in $H^1_\ast(\Gamma,F)$ by Proposition 8 of \cite{Serre1}, which is essentially a consequence of Chebotarev's Density
Theorem. Because $H^1(k,-)$ commutes with limits of finite Galois modules and by left exactness of the limit functor we can deduce that $\Scha^1(k,\Tell M)$ is contained in $H^1_\ast(\Gamma,\Tell M)$, and $H^1_\ast(\Gamma,\Tell M)$ is isomorphic to $H^1_\ast(L^M,\Tell M)$ by \cite{Serre1} Proposition 6, hence the canonical injection. By \emph{loc.cit.}\ Proposition 9 the $\IZ_\ell$--module $H^1(L^M,\Tell M)$ is finitely generated, and we have an isomorphism of finite dimensional vector spaces 
$$H^1(L^M,\Tell M) \otimes \IQ_\ell \cong H^1(L^M,\Vell M)$$
This identifies $H^1_\ast(L^M,\Tell M)\otimes\IQ_\ell$ with a subspace of $H^1_\ast(L^M,\Vell M)$, which in turn is a subspace of $H^1_\ast(\fl^M,\Vell M)$ by Lemma \ref{Lem:LazardH1ast}.
\end{proof}

\vspace{4mm}
\begin{lem}\label{Lem:LinAlgBaseChange}
Let $K_1|K_0$ be an extension of fields (think of $\IQ_\ell|\IQ$). Let $E_0$ be a $K_0$--algebra, let $D_0$ and $V_0$ be $E_0$--modules and let $X_0$ be a $K_0$--linear subspace of $D_0$. Denote by $E_1$, $D_1$, $V_1$ and $X_1$ the corresponding objects
obtained by tensoring with $K_1$. Define
\begin{eqnarray*}
\overline X_0 & := & \{x\in D_0 \tq f(x) \in f(X_0) \pt f\in \Hom_{E_0}(D_0, V_0)\}\\
\overline X_1 & := & \{x\in D_1 \tq f(x) \in f(X_1) \pt f\in \Hom_{E_1}(D_1, V_1)\}
\end{eqnarray*}
Then, the inclusion $\overline X_1 \subseteq \overline X_0\otimes K_1$ holds. In particular, if the equality $X_0=\overline X_0$ holds, then the equality $X_1=\overline X_1$ holds as well. 
\end{lem}

\begin{proof}
Let $x$ be an element of $\overline X_1 \subseteq D_1$ and let us show that $x$ belongs to $\overline X_0\otimes K_1$. Every $E_0$--linear map $D_0\to V_0$ gives rise by $K_1$--linear extension to a $E_1$--linear map $D_1\to V_1$, so by definition of $\overline X_1$ there exists in particular for every $f \in \Hom_{E_0}(D_0, V_0)$ an element $x^f\in X_1$ such that $f(x)=f(x^f)$. Let $(t_i)_{i\in I}$ be a $K_0$--basis of $K_1$, so we can write $x$ and $x^f$ as a sums
$$x = \sum_{i\in I}x_i\otimes t_i \qqet x^f = \sum_{i\in I}x_i^f\otimes t_i$$
for unique elements $x_i\in D_0$ and $x_i^f \in X_0$, almost all zero. We have to show that the $x_i$ belong to $\overline X_0$ for all $i\in I$. The equality $f(x)=f(x^f)$ reads 
$$\sum_{i\in I}f(x_i)\otimes t_i = \sum_{i\in I}f(x_i^f)\otimes t_i$$
Linear independence of the $t_i$'s over $K_0$ implies that we have in fact $f(x_i)=f(x_i^f)$ for all $i$. Hence for every $i\in I$ and every $f\in \Hom_{E_0}(D_0, V_0)$ we have $f(x_i)\in f(X_0)$, that is, $x_i \in \overline X_0$ as we wanted to show. As for the additional statement, if we have $X_0=\overline X_0$ then the inclusions 
$$X_0 \otimes K_1 \:\:\overset{\mathrm{def}}{=}\:\: X_1 \:\:\subseteq\:\: \overline X_1 \:\:\subseteq\:\: \overline X_0 \otimes K_1$$
must all be equalities.
\end{proof}

\vspace{4mm}
\begin{lem}\label{Lem:LinearAlgebra1}
Let $K$ be a field of characteristic zero, let $E$ be a finite product of finite dimensional division algebras over $K$, let $D$ and $V$ be finite dimensional $E$--modules, and suppose that $V$ is faithful. Let $X$ be a $K$--linear subspace of $D$. An element $v\in D$ belongs to $X$ if and only if $f(v)$ belongs to $f(X)$ for all $E$--linear maps $f:D\to V$.
\end{lem}

\begin{proof}
We only show the case where $E$ is a division algebra over $K$, the proof of the general case is similar. That $V$ is faithful means then just that $V$ is nonzero, and without loss of generality we may suppose that $V$ is $E$, so we are considering $E$--linear forms $f: D\to E$. Let $\tr_{E|K}:E\to K$ be a trace map, which for our purpose can be just any $K$--linear map with the property 
$$\tr_{E|K}(yx)=0 \pt y\in E \implies x=0$$
Such a trace map always exists (see e.g. \cite{GiSza06} section 2.6). Consider then the $K$--linear map
\begin{eqnarray*}
\Hom_E(D,E) & \to & \Hom_K(D,K)\\
f & \mapsto & \tr_{E|K}\circ f
\end{eqnarray*}
We claim that this is an isomorphism of $K$--vector spaces. We only have to show injectivity, surjectivity follows then by dimension counting. To show injectivity, we can suppose that $D = E$. The above map sends then an $E$--linear endomorphism of $E$, which is just multiplication on the right by some $x\in E$ to the $K$--linear map $y\mapsto \tr(yx)$. If this map is zero, then $x$ must be zero by the above property of the trace, hence injectivity. The hypothesis on $v$ implies that
$$\tr_{E|K}f(v) \in \tr_{E|K}f(X)$$
for all $f\in\Hom_E(D,E)$, hence $f(v) \in f(X)$ for all $f\in \Hom_K(D,K)$, hence $v\in X$ by standard linear algebra.
\end{proof}

\vspace{4mm}
\begin{proof}[Proof of Theorem \ref{Thm:FinitenessEssentials}]
By Proposition \ref{Pro:SchaAndH1astRank}, it is enough to show that the vector space $H^1_\ast(\fl^M,\Vell M)$ is trivial. Set $E_\ell := \End_{\overline k}(G) \otimes\IQ_\ell$ and $X_\ell := \im(u)\otimes\IQ_\ell$, denote by $D_\ell$ the $E_\ell$--submodule of $G(\overline k)\otimes\IQ_\ell$ generated by $X_\ell$, and define
$$\overline X_\ell := \big\{x\in D_\ell \:\big|\: f(x)\in f(X_\ell) \pt f\in\Hom_{E_\ell}(D_\ell, \Vell G) \big\}$$
By Theorem \ref{Thm:LieCohoMain} we have to check that the equality $X_\ell = \overline X_\ell$ holds. Fix an embedding of $k$ into the field of complex numbers $\IC$. Set $\Vnot G := H_1(G(\IC), \IQ)$ and $E_0 := \End_{\overline k}(G)\otimes\IQ$ and $X_0 := \im(u)\otimes\IQ$, and denote by $D_0$ the $E_0$--submodule of $G(\overline k)\otimes\IQ$ generated by $X_0$. Note that $\Vnot G$ is a faithful $E_0$--module, and that there is a natural isomorphism $\Vell G \cong \Vnot G\otimes\IQ_\ell$. By Lemma \ref{Lem:LinAlgBaseChange} it is now enough to check the equality $X_0 = \overline X_0$ for
$$\overline X_0 := \big\{x\in D_0 \:\big|\: f(x)\in f(X_0) \pt f\in\Hom_{E_0}(D_0, \Vnot G) \big\}$$
By hypothesis, the $\IQ$--algebra $E_0$ is a product of division algebras, hence the equality $X_0 = \overline X_0$ indeed holds by Lemma \ref{Lem:LinearAlgebra1}.
\end{proof}

\vspace{4mm}
\begin{para}\label{Par:MoreLinearAlgebra}
One can think of other linear algebra conditions on the objects $E$, $D$, $V$ and $X$ than those in Lemma \ref{Lem:LinearAlgebra1} which ensure the equality $X = \overline X$. For instance, the conclusion of the Lemma holds true for any finite dimensional semisimple algebra $E$ over $K$ and faithful $V$, if $X$ is of dimension $\leq 1$, or if $X$ is an $E$--submodule or $D$. One can conclude along the same lines that if $M = [u:Y \to G]$ is a 1--motive where $G$ is an abelian variety or a torus, such that the image of $u$ generates an $\End_{\overline k} \otimes\IQ$--submodule of $G(\overline k)\otimes\IQ$ or such that $u(Y)$ is of rank $\leq 1$, then $\Scha^1(k,\Tell M)$ is finite.
\end{para}

\vspace{4mm}
\begin{para}
Our strategy of proving finiteness of $\Scha^1(k,\Tell M)$ consisted of showing that the a priori larger group $H^1_\ast(L^M, \Tell M)$ is finite. This strategy does not succeed always, indeed, there exist 1--motives $M$ such that the group $H^1_\ast(L^M, \Tell M)$ is infinite, yet $\Scha^1(k,\Tell M)$ is finite. The point here is that $H^1_\ast(\Gamma, \Tell M)$ only sees the primes at which $\Tell M$ is unramified, whereas $\Scha^1(k,\Tell M)$ sees all primes. 
\end{para}

\vspace{14mm}
\section{\texorpdfstring{The torsion of $\Scha^1(k,\Tell M)$}{The torsion of Scha1(k, Tl M)}} \label{Sec:TorsionOfSchaTellM}

\begin{par}
In this section we complete the proof of Theorem \ref{Thm:Main3} by examining the finite torsion part of the group $\Scha^1(k,\Tell M)$. The key ingredient for this is the following abstract lemma. 
\end{par}

\vspace{4mm}
\begin{lem}\label{Lem:KeyLem}
Let $T$ be a finitely generated free $\IZ_\ell$--module and set $V := T \otimes\IQ_\ell$. Let $D\subseteq L \subseteq \GL(T)$ be Lie subgroups with Lie algebras $\fd$ and $\fl$ respectively. If 
\begin{enumerate} 
 \item the set $\{\pi\circ x \tq x\in\fl, \pi\in V^\ast\}$ is a linear subspace of $V^\ast$, and
 \item for all open subgroups $H\subseteq L$ containing $D$ the equality $T^H=T^L$ holds,
\end{enumerate}
then the map $r:H^1_\ast(L,T)\to H^1_\ast(D,T)$ given by restriction of cocycles is injective on torsion elements.
\end{lem}

\begin{par}
This generalises Lemma 4.1 in \cite{Jos10}, which we get back by taking for $D$ the trivial group. In our application, $T$ will be $\Tell M$ for a 1--motive $M$, $L$ will be $L^M$, i.e. the image of $\Gamma := \Gal(\overline k|k)$ in $\GL(\Tell M)$, and $D$ will be the image in $\GL(\Tell M)$ of a decomposition group $D_v \subseteq \Gamma$. 
\end{par}

\begin{proof}[Proof of Lemma \ref{Lem:KeyLem}]
\begin{par}
Let $c:L\to T$ be a cocycle representing an element of order $\ell$ in $\ker(r)$, and let us show that $c$ is a coboundary. Because $c$ represents a torsion element in $H^1(L,T)$ its image in $H^1(L,V)$ is trivial. Thus, identifying $T$ with a subset of $V$, there exists $v\in V$ with the property that $c(g)=gv-v$ for all $g\in L$. The cocycle $c$ is a coboundary if $v$ belongs to $v\in V^L+T$, and that is what we will show. 
\end{par}

\begin{par} {\bf Claim.}
\emph{We claim that $v$ belongs to $(T+V^D) \cap (T + V^\fl)$.}
\end{par}

\begin{par}                                                                                                                    
Since the restriction of $c$ to $D$ is a coboundary, there exists $t\in T$ such that $c(g) = gt-t$ for all $g\in D$, hence $v-t \in V^D$ and we have indeed $v\in T+V^D$. To say that the cohomology class of $c$ belongs to $H^1_\ast(L,T)$ is to say that for each $g\in L$ there exists an element $t_g\in T$ such that $c(g) = gt_g-t_g$.  Let $N$ be an open normal subgroup of $\fl$ on which the exponential map $\exp: N\to \fl$ is defined, so that $V^{\angl g} = \ker(\exp(g))$ for all $g\in N$. We then have
$$v \in\:\: \bigcap_{g\in L}(T + V^{\angl g}) \:\:\subseteq\:\:\bigcap_{g\in N}(T + V^{\angl g}) \:\:=\:\: \bigcap_{x\in \fl}(T + \ker(x))$$
By  Lemma 4.4 in \cite{Jos10}, which applies because of the hypothesis (1), we have 
$$v \in \:\: T + \bigcap_{x\in \fl}\ker(x) \:\:=\:\: T+V^\fl$$
hence the claim.                                                                                                                                       
\end{par}
\begin{par}
By modifying $v$ by an element of $T$ we may suppose without loss of generality that $v$ belongs to $V^D$, and in particular to $V^\fd$. The finite group $G := D/(N\cap D)$ acts on $V^\fd$ as well as on $V^\fl$. By Maschke's theorem there exists a $\IQ_\ell$--linear, $G$--equivariant retraction map $r:V^\fd\to V^\fl$ of the inclusion $V^\fl\to V^\fd$. Restricting $r$ to $V^\fl + (T\cap V^\fd)$ we find a decomposition of $G$--modules
$$V^\fl + (T\cap V^\fd) = V^\fl \oplus \big(\ker r \cap(T\cap V^\fd)\big)$$
Writing $v=v_1+t_1$ with $v_1\in V^\fl$ and $t_1 \in \ker r \cap T\cap V^\fd$ according to this decomposition we see that $v_1$ (and also $t_1$) is fixed under $G$ because $v$ is so, hence we have
$$v\in (V^\fl \cap (V^\fd)^G)+T = (V^N \cap V^D)+T = V^{ND} + T$$
The subgroup $ND$ of $L$ is open and contains $D$, hence $v\in V^L+T$ by hypothesis (2).
\end{par}
\end{proof}

\vspace{4mm}
\begin{lem}\label{Lem:TrivialScha}
Let $M=[u:Y\to G]$ be a 1--motive where $G$ is an abelian variety or a torus, such that $\End_{\overline k}(G)\otimes\IQ$ is a product of division algebras. If the Galois action on $Y$ is trivial, then $\Scha^1(k,\Tell M)$ is trivial.
\end{lem}

\begin{proof}
For every finite Galois extension $k'|k$ we have $H^0(k,\Tell M)=H^0(k',\Tell M) \cong H^{-1}(M)\otimes\IZ_\ell$ by Corollary \ref{Cor:TellFixGlobal}. Hence we have $(\Tell M)^{L^M}=(\Tell M)^U$ for all open subgroups $U$ of $L^M$. It follows from Proposition 3.1 and 3.2 of \cite{Jos10}, which uses the hypothesis on $\End_{\overline k}(G)\otimes\IQ$ that the image of the bilinear map
$$\fl^M\times (\Vell M)^\ast \to (\Vell M)^\ast \qquad \qquad (x,\pi) \mapsto \pi\circ x$$
is a linear subspace of $(\Vell M)^\ast$. The hypothesises of Lemma \ref{Lem:KeyLem} are thus satisfied, and taking for $D$ the trivial group, it shows that $H^1_\ast(L^M, \Tell M)$ is torsion free. By Theorem \ref{Thm:FinitenessEssentials}, this group is also finite, hence trivial, and we conclude by Proposition \ref{Pro:SchaAndH1astRank}.
\end{proof}

\vspace{4mm}
\begin{proof}[Proof of the Theorem 3]
Let $M=[u:Y\to G]$ be a 1--motive over $k$. We have constructed the pairing of the Main Theorem, and shown in Corollary \ref{Cor:Scha0Scha2PairingDef} that it is nondegenerate on the left, and in \ref{Pro:CalculScha0} that $\Scha^0(k,M)$ is finite. Suppose then that $G$ is an abelian variety or a torus, such that $\End_{\overline k}(G)\otimes\IQ$ is  a product of division algebras. By Corollary \ref{Cor:Scha0Scha2PairingDef}, it remains to prove that the canonical map 
\begin{equation}
\Scha^0(k,M)\otimes \IZ_\ell \to \Scha^1(k,\Tell M)\tag{$\ast$}
\end{equation}
constructed in \ref{Par:CanonicalAst} is an isomorphism. We define $Z :=H^{-1}(M)=\ker u$ and use Proposition \ref{Pro:CalculScha0} to identify $\Scha^1(k,M)\otimes \IZ_\ell$ with $\Scha^1(k, Z \otimes \IZ_\ell)$. Fix a finite Galois extension $k'|k$ over which $Z$ is constant. For every place $w$ of $k'$ we write $k_w$ for the completion of $k$ at the restriction of $w$ to $k$, and $D_w := \Gal(k'_w|k_w)$. From the Hochschild--Serre spectral sequence we get a commutative diagram with exact rows 
$$\begin{diagram} \setlength{\dgARROWLENGTH}{4mm} 
\node{0}\arrow{e} \node{H^1(\Gal(k'|k), Z \otimes\IZ_\ell)} \arrow{s}\arrow{e}\node{H^1(k,\Tell M)} \arrow{s}\arrow{e} \node{H^1(k',\Tell M)} \arrow{s}\\
\node{0}\arrow{e} \rnode{\prod_w H^1(G_w, H^0(k_w', \Tell M))} \arrow{e} \rnode{\prod_w H^1(k_w,\Tell M)} \arrow{e} \rnode{\prod_w H^1(k'_w, \Tell M)} 
\end{diagram}$$
The right hand vertical map is injective by Lemma \ref{Lem:TrivialScha}, hence every element of $\Scha^1(k,\Tell M)$
comes from a unique element of $H^1(\Gal(k'|k), Z \otimes\IZ_\ell)$, hence from $H^1(k, Z\otimes\IZ_\ell)$. It remains to show that this element is in $\Scha^1(k,Z\otimes \IZ_\ell)$. To this end, we consider the following diagram
$$\begin{diagram}
\setlength{\dgARROWLENGTH}{4mm}
\node{0}\arrow{e} \node{H^1(k, Z \otimes\IZ_\ell)} \arrow{se,t,--}{\delta}\arrow{s}\arrow{e}\node{H^1(k,\Tell M)}\arrow{s}\\
\node[2]{\prod H^1(k_v,Z \otimes\IZ_\ell)} \arrow{e} \node{\prod H^1(k_v,\Tell M)}
\end{diagram}$$
where the horizontal maps are induced by the morphism of 1--motives $[Z\to 0] \to [Y\to G]$. Injectivity of the top horizontal map follows from Corollary \ref{Cor:TellFixGlobal}. We have thus $\ker\delta \cong \Scha^1(k,\Tell M)$ and must show that every element of $\ker\delta$ maps already to zero in $H^1(k_v,Z \otimes\IZ_\ell)$ for all $v\in \Omega$, that is, $\ker\delta = \Scha^1(k,Z\otimes\IZ_\ell)$.\\
Fix an element $x$ of $\Scha^1(k,\Tell M)$, a place $v$ and let $D_v$ be a decomposition group for $v$. We know that $x$ comes via inflation from an element $z$ of the finite group $H^1_\ast(L^M, \Tell M)$. Write $D$ for the image of $D_v$ in $\GL(\Tell M)$. This $D$ is a Lie subgroup of $L^M$, and by hypothesis $z$ restricts to zero in $H^1(D, \Tell M)$. By Lemma
\ref{Lem:KeyLem} (using again Proposition 3.1 and 3.2 of \cite{Jos10}) we conclude that there is an open subgroup $U$ of $L^M$ containing $D$, such that $z$ is already zero in $H^1(U, \Tell M)$. This shows as well that there is an open subgroup $\Gamma'$ of $\Gamma$ containing $D_v$ such that $x$ maps to zero in $H^1(\Gamma', \Tell M)$. Consider then the diagram
$$\begin{diagram}
\setlength{\dgARROWLENGTH}{4mm}
\node{0}\arrow{e} \node{H^1(\Gamma, Z \otimes\IZ_\ell)} \arrow{se,t,--}{\delta'}\arrow{s}\arrow{e}\node{H^1(\Gamma,\Tell M)}\arrow{s,r}{x\mapsto 0}\\
\node{0}\arrow{e} \node{H^1(\Gamma', Z \otimes\IZ_\ell)} \arrow{s}\arrow{e}\node{H^1(\Gamma',\Tell M)}\arrow{s}\\
\node[2]{H^1(D_v,Z \otimes\IZ_\ell)} \arrow{e} \node{H^1(D_v,\Tell M)}
\end{diagram}$$
We know that the element $x\in \Scha^1(k,\Tell M)$ comes from an element of $\ker\delta'$. The middle row is exact by Corollary  \ref{Cor:TellFixGlobal} and because $\Gamma'$ is the Galois group of a number field, so that this element maps to zero in $H^1(\Gamma', Z \otimes\IZ_\ell)$, hence in $H^1(D_v,Z \otimes\IZ_\ell)$.
\end{proof}

\vspace{14mm}
\section{Tate 1--motives and Leopoldt's conjecture}\label{Sec:Leopoldt}

\begin{par}
In this section we study the pairing of Theorem \ref{Thm:Main3} in the case where $M$ is a \emph{Tate 1--motive} over $k$, that is, a 1--motive of the form $M=[\IZ^r \to \IG_m^s]$. For such $1$--motives, and more generally for 1--motives of the form $[Y\to G]$ where $G$ is a torus, I suspect that the pairing of the Main Theorem is perfect. For a Tate 1--motive this amounts to say that $\Scha^2(k,M)$ is trivial. I will show the following sharper version of theorem \ref{Thm:Main2} stated in the itroduction.
\end{par}

\vspace{4mm}
\begin{thm}\label{Thm:DualityAndLeopoldt} 
Let $k$ be a number field with ring of integers $\cO_k$, and let $\ell$ be a prime number. If for every $1$--motive of the form $M=[\IZ^r \to \IG_m^2]$ over $\spec(\cO_k)$ the group $\Scha^2(k,M^\vee)[\ell^\infty]$ is trivial, then the statement of Leopoldt's conjecture is true for $k$ and $\ell$.
\end{thm}

\vspace{4mm}
\begin{para}
We will work with the following formulation of Leopoldt's conjecture (\cite{Neukirch00}, statement (iii) of Theorem 10.3.6). For a finite prime $\fp$ of $k$, let $\cO_{k,\fp}$ denote the ring of integers of the completion of $k$ at $\fp$. There is a canonical map 
$$i_\ell: \cO_k^\ast \otimes_\IZ \IZ_\ell \to \prod_{\fp|\ell}
\cO_{k,\fp}^\ast\hotimes \IZ_\ell$$
which on each component $i_\fp: \cO_k^\ast\otimes \IZ_\ell \to
\cO_{k,\fp}\hotimes \IZ_\ell$ is obtained by applying $-\hotimes \IZ_\ell$ to the inclusion $\cO_k^\ast \subseteq \cO_{k,\fp}^\ast$. Leopoldt's conjecture asserts that the map $i_\ell$ is injective. Note that $i_\ell$ is injective on torsion elements, and injective if $\cO_k^\ast$ is of rank $\leq 1$. 
\end{para}

\vspace{4mm}
\begin{proof}[Proof of theorem \ref{Thm:DualityAndLeopoldt}]

\begin{par}
We suppose Leopoldt's conjecture is false for $k$ and $\ell$, so there exists a nontorsion element $z\in \ker(i_\ell) \subseteq \cO_k^\ast \otimes_\IZ\IZ_\ell$, which we may write as
$$z = \sum_{i=1}^n \eps_i \otimes \lambda_i$$
where $n \geq 2$ is the rank of $\cO_k^\ast$ and $\eps_1, \ldots, \eps_n$ are $\IZ$--linearly independent elements of $\cO_k^\ast$. By reordering the $\eps_i$'s and replacing $\eps_1$ by $\eps_1^{-1}$ if necessary, we may as well assume $\lambda_1 + \lambda_2 \neq 0$. We will now construct a $1$--motive $M$ of the form $M = [u:\IZ^{2n-1}\to \IG_m^2]$ over $\spec(\cO_k)$ such that the group $\Scha^1(k,\Tell M)$ is infinite. The 1--motive dual to $M$ is then of the form $M^\vee = [\IZ^2\to\IG_m^{2n-1}]$ and $\Scha^2(k,M^\vee)$ will be infinite by Proposition \ref{Pro:Scha1TateModScha2}. Let $Y \simeq \IZ^{2n-1}$ be the group matrices 
$$y = \bigg(\!\!\!\begin{array}{ccccc} y_{11} &y_{12} &y_{13} &\cdots& y_{1n}\\y_{21} &y_{22} &y_{23} &\cdots & y_{2n}\end{array}\!\!\!\bigg)$$
with integer coefficient satisfying $y_{11}+y_{22}=0$, and define the morphism $u$ by
$$u(y)= \bigg(\!\!\!\begin{array}{c} y_{11}\eps_1 + y_{12}\eps_2 + y_{13}\eps_3
+ \cdots + y_{1n}\eps_n\\ y_{21}\eps_1 + y_{22}\eps_2 + y_{23}\eps_3 +  \cdots +
y_{2n}\eps_n \end{array}\!\!\!\bigg) \quad \in \IG_m^2(\cO_k)$$
where we decided to write the group $\IG_m(\cO_k) = \cO_k^\ast$ additively. So if $\eps$ denotes the column vector of the $\eps_i$'s, we have just $u(y) = y\eps$. We will prove the following lemma later:
\end{par}

\begin{lem}\label{Lem:SchinzelLeopoldt}
For each $1\leq i \leq n$, there exists $y\in Y$ such that $\binom{\eps_i}{\eps_i} \equiv u(y)\bmod\fp$ holds in $\IG_m^2(\kappa_\fp)$, where $\kappa_\fp = $ is the residue field at $\fp$. 
\end{lem}

\begin{par}
Set $U := \spec(\cO_k[\ell^{-1}])$ and denote by $c_i$ and $c$ respectively the image of $\binom{\eps_i}{\eps_i} \otimes 1$ and $\binom zz$ under the composite map 
$$\IG_m^2(U)\otimes \IZ_\ell \to H^0(U,M) \otimes \IZ_\ell \to H^1(U,\Tell M)$$
where the first map is induced by the projection $\IG_m^2(U) \to \IG_m^2(U)/u(Y)\cong H^0(U,M)$ and the second map is the injection defined in \ref{Pro:FiltrationTate}. The $\IZ_\ell$--submodules of $\IG_m^2(\cO_k)\otimes\IZ_\ell$ generated by $\binom zz$ and by $u(Y)$ intersect trivially because $\lambda_1+ \lambda_2 \neq 0$, hence $c$ is of infinite order in $H^1(U,\Tell M) \subseteq H^1(k,\Tell M)$. I claim that $c$ belongs to $\Scha^1(k,\Tell M)$. Fix a place $\fp$ of $k$ of residual characteristic $p$, and let us show that the restriction of $c$ to $H^1(k_\fp,\Tell M)$ is zero. In the case $p = \ell$ this is true by construction, considering the commutative diagram
$$\begin{diagram}
\setlength{\dgARROWLENGTH}{4mm}
\node{(\cO_k^\ast)^2 \otimes \IZ_\ell}\arrow{s}\arrow{e}\node{H^1(U,\Tell M)}\arrow{s}\\  
\node{(\cO_{k,\fp}^\ast)^2 \hotimes \IZ_\ell}\arrow{e}\node{H^1(k_\fp, \Tell M)}
\end{diagram}$$
and that the image of $\binom zz$ is already zero in $(\cO_{k,\fp}^\ast)^2\hotimes\IZ_\ell$. Suppose now that $p\neq \ell$, so $\Tell M$ is unramified at $\fp$. Because $c= \lambda_1c_1 + \cdots + \lambda_nc_n$, it suffices to show that the restriction of each $c_i$ to $H^1(\kappa_\fp,\Tell M)$ is zero. In view of the commutative diagram
$$\begin{diagram}
\setlength{\dgARROWLENGTH}{4mm}
\node[2]{(\cO_k^\ast)^2 \otimes \IZ_\ell}\arrow{s}\arrow{e}\node{H^1(U,\Tell M)}\arrow{s}\\  
\node{Y\otimes \IZ_\ell}\arrow{e,t}{u\bmod \fp}\node{(\kappa_\fp^\ast)^2 \otimes \IZ_\ell}\arrow{e}\node{H^1(\kappa_\fp, \Tell M)}
\end{diagram}$$
this amounts to show that there exists $y\in Y$ such that $\binom{\eps_i}{\eps_i} \equiv u(y)\bmod\fp$ holds in $\IG_m^2(\kappa_\fp)$, which is what we claimed in Lemma \ref{Lem:SchinzelLeopoldt}. Hence $c$ belongs indeed to $\Scha^1(k,\Tell M)$ and is of infinite order, and thus $\Scha^1(k,\Tell M)$ is infinite.
\end{par}
\end{proof}

\vspace{4mm}
\begin{proof}[Proof of Lemma \ref{Lem:SchinzelLeopoldt}]
Fix $1\leq i \leq n$ and a maximal ideal $\fp$ of $\cO_k$ with residue field $\kappa_\fp$. We have to find a matrix $y \in Y$ such that $\binom{\eps_i}{\eps_i}$ is congruent to $u(y)$ modulo $\fp$. For $i\neq 1,2$ such a $y$ exists trivially. Let $J_1 \subseteq \IZ$ be the ideal consisting of those $m\in \IZ$ such that $m\eps_1 \bmod \fp$ is in the subgroup of $\kappa_\fp^\ast$ generated by $\eps_2$ and let $a_1\geq 1$ be the positive generator of $J_1$. Similarly, define $J_2$ and $a_2$. There exists $b_1, b_2 \in \IZ$ such that the linear dependence relations
$$a_1\eps_1 + b_2\eps_2 = 1 \qqet b_1\eps_1 + a_2\eps_2 = 1$$
hold in the finite group $\kappa_\fp^\ast$, written additively. Note that $b_i$ is a multiple of $a_i$. We claim that the integers $a_1$ and $a_2$ are coprime. Indeed, suppose there exists a prime $\ell$ dividing $a_1$ and $a_2$ so that we can write $a_i = \ell a_i'$ and $b_i = \ell b_i'$. Let $Z$ be the subgroup of $\kappa_\ell^\ast$ generated by $\eps_1$ and $\eps_2$. Since 
$\kappa_\fp^\ast[\ell]$ is cyclic of order $\ell$, we may suppose that $\kappa_\fp^\ast[\ell] \cap Z$ is contained in the subgroup of $\kappa_\ell^\ast$ generated by, say, $\eps_2$. Thus, the point 
$$T := a_1'\eps_1 + b_2'\eps_2 \qquad \in \kappa_\fp^\ast[\ell] \cap Z$$
can be written as $T = c\eps_2$, and we get the relation $a_1'\eps_1 + (b_2'-c)\eps_2 = 1$ which contradicts the minimality of $a_1$. Therefore, $a_1$ and $a_2$ are coprime as claimed, and we can choose integers $c_1,c_2$ such that $a_1c_1+a_2c_2 = 1$. The matrices
$$y_1 = \bigg(\!\!\!\begin{array}{ccccc} 
1 - a_1c_1 & - c_1b_2 & 0 &\cdots& 0 \\
1 - c_2b_1 & - a_2c_2 & 0 &\cdots & 0\end{array}\!\!\!\bigg) \qqet 
y_2 = \bigg(\!\!\!\begin{array}{ccccc} 
- a_1c_1 & 1 - c_1b_2 & 0 &\cdots& 0 \\
- c_2b_1 & 1 - a_2c_2 & 0 &\cdots & 0\end{array}\!\!\!\bigg)$$
belong to $Y$, and we have $\binom{\eps_i}{\eps_i} \equiv u(y_i)\bmod\fp$ as desired.
\end{proof}

\vspace{14mm}
\section{\texorpdfstring{A semiabelian variety with infinite $\Scha^2$}{A semiabelian variety with infinite Scha 2}}\label{Sec:InfiniteScha2}

\vspace{4mm}
\begin{par}
In this section we prove Theorem \ref{Thm:Main1} by producing a semiabelian variety $G$ over $\IQ$ such that $\Scha^2(\IQ,G)$ contains $\QZ$ as a subgroup, and hence in particular is infinite. The technique is similar to that in the previous paragraph, and here we exploit now that for elliptic curves of sufficiently big rank the statement analogue to Leopoldt's conjecture trivially fails.
\end{par}

\begin{para}
Let $E$ be an elliptic curve over $\IQ$ of rank at least $3$ and let $P_1,P_2,P_3 \in E(\IQ)$ be $\IZ$--linearly independent rational points. Let us write $A$ for the abelian threefold $E^3$ over $\IQ$ and $Y$ for the group of $3\times 3$ matrices of trace zero with integer coefficients. Looking at $Y\simeq \IZ^8$ as a Galois module with trivial Galois action, we consider the 1--motive 
$$M = [u:Y\to A] \qquad \qquad u(y) = yP = \left(\begin{array}{c}
y_{11}P_1+y_{12}P_2 + y_{13}P_3\\
y_{21}P_1+y_{22}P_2 + y_{23}P_3\\
y_{31}P_1+y_{32}P_2 + y_{33}P_3
\end{array}\right) \in E(\IQ)^3 = A(\IQ)$$ 
The map $u$ is injective, and I will use $X$ as a shorthand for the group $u(Y)\subseteq A(\IQ)$. This 1--motive $M$ is of special interest because it produces a counterexample to the so called \emph{problem of detecting linear dependence}: Although $P\notin X$, and even $nP\notin X$ for all $n\neq 0$, there exists for every prime $p$ where $E$ has good reduction an element $x\in X$ such that $P$ is congruent to $x$ modulo $p$. The verification of this is similar to the proof of Lemma \ref{Lem:SchinzelLeopoldt}, see \cite{JPer10}. Using Theorem \ref{Thm:LieCohoMain} one shows that $H^1_\ast(\fl^M,\Vell M)$ is nontrivial --- this is what makes the counterexample work, and also how it was found in the first place. 
\end{para}

\vspace{4mm}
\begin{para}\label{Par:nPdivisibleClaims}
I claim that the Tate--Shafarevich group in degree 2 of the semiabelian variety dual to the 1--motive $M$ constructed in the previous paragraph contains a subgroup isomorphic to $\QZ$. By Proposition \ref{Pro:Scha1TateModScha2} this amounts to say that for each prime number $\ell$ the Tate--Shafarevich group 
$$\Scha^1(\IQ,\Tell M)$$
is of rang $\geq 1$ as a $\IZ_\ell$--module. Fix a prime $\ell$ and let us denote by $[c_P]$ the cohomology class of $P\otimes 1$ via the injection $H^0(\IQ,M)\otimes \IZ_\ell \to H^1(\IQ,\Tell M)$ from proposition \ref{Pro:FiltrationTate}. A cocycle $c_P$ representing $[c_P]$ is explicitly given by
$$c_P(\sigma) = (\sigma P_i-P_i)_{i=0}^\infty$$
where $(P_i)_{i=0}^\infty$ are elements of $A(\overline \IQ)$ such that $P_0=P$ and $\ell P_i = P_{i-1}$. Up to a coboundary, $c_P$ does not depend on the choice of the division points $P_i$. As the class $[P]$ of $P$ in $H^0(\IQ,M) \cong A(\IQ)/X$ is of infinite order, the element $[c_P]\in H^1(\IQ,\Tell M)$ is of infinite order too. We claim that $n[c_P]$ belongs to $\Scha^1(\IQ,\Tell M)$ for some integer $n\geq 1$ (depending on $\ell$). To check this, we must show that for every finite place $p$ of $\IQ$ the restriction of $nc_P$ to a decomposition group $D_p$ is a coboundary. In the case where $\ell=2$ and $p=\infty$ we should also demand that the restriction of $n[c_P]$ to $H^1(\Gal(\IC|\IR),\Tell M)$ is zero, but we can ignore this by choosing $n$ to be even. So from now on we will stick to finite primes $p$ only.
\end{para}

\vspace{4mm}
\begin{lem}\label{Lem:CoboundaryIffDivisible}
Let $p$ be a prime and let $D_p \subseteq \Gal(\overline \IQ|\IQ)$ be a decomposition group at $p$. The restriction of $c_P$ to $D_p$ is a coboundary if and only if the class of $P$ in $A(\IQ_p)/X = H^0(\IQ_p,M)$ is $\ell$--divisible.
\end{lem}

\begin{proof}
Choose an algebraic closure $\overline \IQ_p$ of $\IQ_p$ and an embedding of $\overline \IQ$ into $\overline \IQ_p$ in such a way that the given decomposition group $D_p$ equals $\Gal(\overline \IQ|(\overline \IQ\cap \IQ_p))$. Consider the commutative diagram with exact rows
$$\begin{diagram}
\setlength{\dgARROWLENGTH}{4mm}
\node{0}\arrow{e}\node{H^0(\IQ,M)\otimes\IZ_\ell}\arrow{s}\arrow{e}\node{H^1(\IQ
,\Tell M)}\arrow{s}\\
\node{0}\arrow{e}\node{H^0(\IQ_p,M)\hotimes\IZ_\ell}\arrow{e}\node{H^1(\IQ_p,
\Tell M)}
\end{diagram}$$
The restriction of $c_P$ to $D_p$ is a coboundary if and only if $[P]\otimes 1 \in H^0(\IQ,M)\otimes\IZ_\ell$ maps to zero in $H^0(\IQ_p,M)\hotimes\IZ_\ell$, that is, if and only if the class of $P$ in $H^0(\IQ_p,M)$ is $\ell$--divisible.
\end{proof}

\vspace{4mm}
\begin{lem}\label{Lem:ClosureFiniteIndex}
For every prime $p$, the closure of $X$ in $A(\IQ_p)$ for the $p$--adic topology is an open subgroup of $A(\IQ_p)$ of finite index.
\end{lem}

\begin{proof}
Because $E(\IQ_p)$ has the structure of a compact $p$--adic Lie group of dimension $1$ there exists an open subgroup of $E(\IQ_p)$ isomorphic to $\IZ_p$, and because $E(\IQ_p)$ is compact, any such subgroup has finite index (\cite{Silv86}, Proposition 6.3). We find thus a short exact sequence of profinite groups 
$$0\to \IZ_p^3 \to A(\IQ_p) \to F \to 0$$
for some finite group $F$. Let $m\geq 1$ be an integer annihilating $F$, so that $mX$ is contained in $\IZ^3_p$. The elements 
$$\left(\begin{array}{c} 0\\ mP_1\\ 0 \end{array} \right) \quad 
\left(\begin{array}{c} 0\\ 0\\ mP_2 \end{array} \right) \quad 
\left(\begin{array}{c} mP_3\\ 0\\ 0 \end{array} \right)$$
of $mX \subseteq \IZ_p^3 \subseteq A(\IQ_p)$ are linearly independent over $\IZ_p$, because each $mP_i \in \IZ_p$ is nonzero. The closure of $mX$ in $\IZ_p^3$ contains the $\IZ_p$--submodule generated by these three points, hence is of finite index in $\IZ_p^3$. We conclude that the closure of $X$ in $A(\IQ)$ has finite index. Every closed subgroup of finite index is also open.
\end{proof}

\vspace{4mm}
\begin{para}\label{Par:ConditionsOnn}
We now come to the proof of the claims made in \ref{Par:nPdivisibleClaims}. First of all, let us choose an integer $n\geq 1$ such that the following conditions are met:
\begin{enumerate}
\item[(0)] If $\ell=2$, then $n$ is even.
\item[(1)] For every prime $p\neq \ell$ where $E$ has bad reduction, the point
$nP$ is $\ell$--divisible in $A(\IQ_p)$.
\item[(2)] For $p=\ell$, the point $nP$ belongs to the closure of $X$ in
$A(\IQ_p)$ for the $p$--adic topology. 
\end{enumerate}
Such an integer $n$ exists. Indeed, start with, say, $n=2$, so condition (0) is satisfied. We have already observed that $A(\IQ_p)$ is an extension of a finite discrete group $F$ by $\IZ_p^3$, so by replacing $n$ by some sufficiently high multiple of $n$ we can assure that $nP$ belongs to the subgroup $\IZ^3_p$ of $A(\IQ_p)$ which is $\ell$--divisible. We do this for all the finitely many primes of bad reduction, so condition (1) is met. As for the last condition, we know that the closure of $X$ in $A(\IQ_p)$ has finite index by Lemma \ref{Lem:ClosureFiniteIndex}, so we again replace $n$ by some sufficiently high
multiple if necessary. In order to show that $n[c_P]$ belongs to $\Scha^1(\IQ,\Tell M)$ it remains to show by Lemma \ref{Lem:CoboundaryIffDivisible} that for each prime $p$ the class of $nP$ in $H^0(\IQ,M)=A(\IQ_p)/X$ is $\ell$--divisible. In other words we must show 

\vspace{2mm}
\begin{par}
{\bf Claim:} \emph{For every $i\geq 0$ there exist elements $Q_i\in A(\IQ_p)$ and $x_i\in X$ such that $\ell^iQ_i +x_i= nP$}
\end{par}

\vspace{2mm}
\begin{par}
\noindent We have already ruled out the case $p=\infty$, and for finite $p$ we will distinguish three cases: First, the case where $p$ is a place of good reduction for $E$ and $p\neq \ell$, second, the case where $p$ is a place of bad reduction and $p\neq \ell$ and finally the case $p=\ell$. All but finitely many primes $p$ fall in the first case. For the finitely many primes that remain, the claim will hold by our particular choice of $n$.
\end{par}
\end{para}

\vspace{4mm}
\begin{par}
\underline{Case 1: good reduction at $p$ and $p\neq \ell$.} In this case, we can consider the surjective reduction map $\red_p:E(\IQ_p)\to E(\IF_p)$. Its kernel is isomorphic to $\IZ_p$, so we get a short exact sequence
$$0\to\IZ_p^3 \to A(\IQ_p) \xrightarrow{\:\:\red_p\:\:} A(\IF_p) \to 0$$
By \cite{JPer10} there exists an element $x\in X$ such that $\red_p(P) = \red_p(x)$ in $A(\IF_p)$. Because $\IZ_p$ is uniquely $\ell$--divisible we can define $Q_i := \ell^{-i}n(P-x)$, and get $\ell^iQ_i+nx=nP$.
 \end{par}

\vspace{4mm}
\begin{par}
\underline{Case 2: bad reduction at $p$ and $p\neq \ell$.} Condition (1) in \ref{Par:ConditionsOnn} ensures that $nP$ is $\ell$--divisible in $A(\IQ_p)$ for bad $p\neq \ell$, so the class of $P$ in $A(\IQ_p)/X$ is $\ell$-divisible as well.
\end{par}

\vspace{4mm}
\begin{par}
\underline{Case 3: $p=\ell$.}  For all $i\geq 0$ the subgroup $p^iA(\IQ_p)$ is open in $A(\IQ_p)$, hence by condition (2) in \ref{Par:ConditionsOnn}, the intersection $X\cap (nP+p^iA(\IQ_p))$ is nonempty. But that means that there exists an element $Q_i\in A(\IQ_p)$ and an element $x_i\in X$ such that $p^iQ_i + x_i=nP$, just as needed.
\end{par}

\vspace{14mm}
\section{Open questions and problems}\label{Sec:OpenQuest}

\begin{par}
I present three open arithmetic questions and an elementary problem in linear algebra which so far have defied all attempts of being solved. The first question is about how far finitely generated subgroups of a Mordell--Weil group are detectable by reduction maps. It is a sharpened version of the problem which in the literature it is named the problem of \emph{detecting linear dependence}. 
\end{par}

\begin{que}
Let $G$ be a semiabelian variety defined over a number field $k$, let $X$ be a finitely generated subgroup of $G(k)$. Denote by $\overline X \subseteq G(k)$ the subgroup of those points $P$, such that for almost all finite primes $\fp$ of $k$ the reduction $P\bmod\fp$ belongs to $X \bmod \fp$ in $G(\kappa_\fp)$. Let $M = [u:Y \to G]$ be a 1--motive where $Y$ is constant and $X=u(Y)$. Is it true that the map 
$$\overline X/X\otimes\IZ_\ell \to H^1_\ast(k,\Tell M)$$
induced by the injection $H^0(k,M) \otimes\IZ_\ell \to H^1(k,\Tell M)$ is an isomorphism?
\end{que}

\begin{par}
A positive answer to this question was given in \cite{Jos10} in the case where $G$ is a geometrically simple abelian variety. In this case we know that $H^1_\ast(k,\Tell M)$ is trivial and get a nice local--global principle for subgroups of Mordell--Weil groups. Apart from a few other isolated examples the question remains open, even in the cases where $G$ is an abelian variety or a torus. The second question is similar in nature, but we impose a stronger local condition.
\end{par}

\begin{que}
Let $G$ be a semiabelian variety defined over a number field $k$, let $X$ be a finitely generated subgroup of $G(k)$, and let $P\in G(k)$ be a rational point. Suppose that for \emph{all} finite primes $\fp$ of $k$ the point $P$ belongs to the closure of $X$ in $G(k_\fp)$ for the $p$--adic topology on $G(k_\fp)$. Does then $P$ belongs to $X$?
\end{que}

\begin{par}
We know the answer to be positive if $G$ is a simple abelian variety, and in some other scattered examples. If we could choose the integer $n$ in \ref{Par:ConditionsOnn} \emph{independently of $\ell$}, the answer to the question would be negative. Thirdly, I would like to ask for a converse to theorem \ref{Thm:DualityAndLeopoldt}:
\end{par}

\begin{que}
Let $k$ be a number field for which the statement of Leopoldt's conjecture holds. Is it true that for all 1--motives $M = [Y \to G]$, where $G$ is a torus, the pairing
$$\Scha^0(k,M) \times \Scha^2(k,M^\vee) \to \QZ$$
of Theorem \ref{Thm:Main3} is a prefect pairing of finite groups?
\end{que}

\begin{par}
At last, motivated by the proof of Theorem \ref{Thm:FinitenessEssentials}, let me state a problem in linear algebra that any first year student can understand:
\end{par}

\begin{para}{\bf Problem.}
Let $K$ be a field, and write $E$ for the $K$--algebra of $n\times n$ matrices with coefficients in $K$. Denote by $V$ and $V_0$ the $E$--modules of $n\times m$ and of $n\times m_0$--matrices respectively. Finally, let $W$ be a $K$--linear subspace of $V$, and define 
$$\overline W := \{v\in V \tq f(v)\in f(W) \pt f\in \Hom_E(V,V_0)\}$$
so $\overline W$ is a linear subspace of $V$ containing $W$. Observe that elements of $\Hom_E(V,V_0)$ are just $m\times m_0$--matrices by Schur's Lemma. The problem is to compute $\overline W$. This means, find an algorithm which takes as an input a $K$--basis of $W$ (this will be some finitely many $n\times m$ matrices), and provides a basis of $\overline W$, or equivalently, provides some finitely many $f_1, \ldots, f_r \in \Hom_E(V,V_0)$ such that 
$$\overline W = \{v\in V \tq f_i(v)\in f_i(W) \mbox{ for } i=1,2,\ldots,r\}$$
Changing scalars from $K$ to a bigger field may result in a smaller dimensional $\overline W$, (i.e. the inclusion of \ref{Lem:LinAlgBaseChange} may be strict). Yet, I don't know of a solution to the problem even in the case where $K$ is algebraically closed.
\end{para}

\vspace{14mm}
\providecommand{\bysame}{\leavevmode\hbox to3em{\hrulefill}\thinspace}
\providecommand{\href}[2]{#2}

$$ $$
\hspace{50mm}
\begin{minipage}[]{100mm}
Peter Jossen\\[1mm]
CNRS, UMR 8628, Math\'ematiques, B\^atiment 425,\\
Universit\'e Paris-Sud, 91450 Orsay cedex - FRANCE\\[2mm]
{\tt peter.jossen@gmail.com} 
\end{minipage}

\end{document}